\documentstyle[epsbox,12pt,a4]{amsart} 

\numberwithin{equation}{section}

\newtheorem{theorem}{Theorem}
\newtheorem{lemma}{Lemma}     
\newtheorem{corollary}{Corollary}
\newtheorem{proposition}{Proposition}

\newtheorem{rem}{Remark}

\newcommand*{\C}{\mathbb{C}}
\newcommand*{\R}{\mathbb{R}}
\newcommand*{\Z}{\mathbb{Z}}
\newcommand*{\Q}{\mathbb{Q}}

\begin{document} 

\title[\null]{
Non-positivity of certain functions associated with analysis on elliptic surfaces
}
\author[\null]{Masatoshi Suzuki}
\date{10 May, 2008}

\subjclass[2000]{11G40, 19F27, 11M41.}

\maketitle


\begin{abstract}
In this paper, we study some basic analytic properties 
of the boundary term of Fesenko's two-dimensional zeta integrals. 
In the case of the rational number field, we show that 
this term is the Laplace transform of certain infinite series consisting of $K$-Bessel functions. 
It is known that the non-positivity property of the fourth log derivative of such series 
is a sufficient condition for the Riemann hypothesis 
of the Hasse-Weil $L$-function attached to an elliptic curve. 
We show that such non-positivity is a necessary condition under some technical assumption.  
\end{abstract}

\section{Introduction}\label{intro}

The main interest of the present paper is the Dirichlet series 
with nonnegative coefficients and its poles. 
Throughout the paper we denote by ${\frak c}$ or $\{c(\nu)\}_{\nu \in A}$ 
a sequence of nonnegative real numbers with a discrete index set $A$. 
Unless we specify the discrete index set $A$, we understand that $A={\Bbb N}$.   
For a nonnegative sequence ${\frak c}= \{c(\nu)\}_{\nu \in A}$, 
we denote by $D_{\frak c}(s)$ the formal Dirichlet series 
\begin{equation}\label{000}
D_{\frak c}(s) = \sum_{\nu \in A} c(\nu) \nu^{-s}. 
\end{equation}
Our basic assumption for $\frak c$ is that $D_{\frak c}(s)$ 
converges (absolutely) on some right-half plane.  
Denote by $\sigma_0 < \infty$ the abscissa of (absolute) convergence of $D_{\frak c}(s)$. 
By well-known Landau's theorem for a Dirichlet series, 
$D_{\frak c}(s)$ has a singularity at $s=\sigma_0$, 
since ${\frak c}$ has a single sign. 
In general, the location of zeros or poles of $D_{\frak c}(s)$ 
in the left of the line $\Re(s)=\sigma_0$ is mysterious and difficult thing
even if it is continued meromorphically to a left of the line $\Re(s)=\sigma_0$. 

In this paper we study one approach to study the poles of $D_{\frak c}(s)$ 
from the viewpoint of the {\it boundary term} 
which is introduced in the theory of 
Fesenko's two-dimensional zeta integrals in~\cite{Fe2}. 
We explain that the so called the {\it single sign property} 
of the theory of boundary term is related deeply to 
the location of the poles of $D_{\frak c}(s)$. 
\smallskip

Let ${\cal E}$ be a two-dimensional arithmetic scheme which is a proper regular model 
of an elliptic curve $E/{\Q}$. 
Let $\zeta_{\cal E}(s)$ be the arithmetic Hasse zeta function of ${\cal E}$,  
which is defined by an Euler product over all closed point of ${\cal E}$. 
Using the computation of the Hasse zeta function for curves over finite fields 
and the description of geometry of models in \cite[Ch.9, Ch.10]{Li}, we obtain 
\begin{equation}\label{001}
\zeta_{\cal E}(s) = 
n_{\cal E}(s)\zeta_E(s), \quad 
\zeta_E(s)=\frac{\zeta(s)\zeta(s-1)}{L(E,s)}
\end{equation}
on $\Re(s)>2$, where $\zeta(s)$ is the Riemann zeta function, 
$L(E,s)$ is the $L$-function of $E/{\Q}$ 
and $n_{\cal E}(s)$ is the product of finitely many, say $J$, Euler factors determined by 
singular fibres of $\cal E$:
\begin{equation} \label{002_1}
n_{\cal E}(s) = \prod_{1 \leq j \leq J}(1-q_j^{1-s})^{-1} 
\quad (\text{$q_{j_1} \not= q_{j_2}$ if $j_1 \not= j_2$}).
\end{equation}
In particular $n_{\cal E}(s)^{\pm 1}$ are holomorphic functions on $\Re(s) > 1$. 
For the two-dimensional arithmetic scheme $\cal E$ and a set $S$ of curves on $\cal E$, 
Fesenko defined a zeta integral of a function on its adelic space 
and a character of its $K_2$-delic group~\cite[\S3]{Fe2}. 
Calculating the zeta integral in two ways, he obtained the following formula for $\Re(s)>2$:
\begin{equation}\label{002}
\widehat{\zeta}(s/2)^2 \cdot c_{\cal E}^{1-s} \cdot \zeta_{\cal E}(s)^2 
= \xi_{\cal E}(s) + \xi_{\cal E}(2-s) + \omega_{\cal E}(s),
\end{equation}
where $\widehat{\zeta}(s)$ is the completed Riemann zeta function $\pi^{-s/2}\Gamma(s/2)\zeta(s)$. 
Here $\xi_{\cal E}(s)$ is an entire function and $c_{\cal E}$ is the constant given by
\[
c_{\cal E} = q_E \prod_{1 \leq j \leq J}q_j
\]
using the conductor $q_E$ of $E$ and values $q_j$ in \eqref{002_1}.
See~\cite[sec.\,40, sec.\,45]{Fe2} for details. 
The third term $\omega_{\cal E}(s)$ in \eqref{002} is called the {\it boundary term}, 
because it can be expressed as an integral 
over the boundary of some two-dimensional adelic object. 
The boundary term $\omega_{\cal E}(s)$ is holomorphic on the right-half plane $\Re(s)>2$. 

Equality $\eqref{001}$ and $\eqref{002}$ 
show that the study of poles of $\omega_{\cal E}(s)$ plays an essential role 
for the study of the zeros of $L(E,s)$. 
The boundary term $\omega_{\cal E}(s)$ has the following integral representation for $\Re(s)>2$
\begin{equation}\label{004}
\omega_{\cal E}(s) 
=  \int_{0}^{1} h_{\cal E}(x) x^{s-2} \frac{dx}{x} 
=  \int_{0}^{\infty} e^{2t} \, h_{\cal E}(e^{-t}) \, e^{-st} dt,
\end{equation}
where $h_{\cal E}(x)$ is a real valued function on $(0,\infty)$. 
Hence the location of poles of $\omega_{\cal E}(s)$ 
is closely related to the behavior of $h_{\cal E}(x)$ as $x$ tends zero.  

For $a,b \in {\Bbb R}_{>0}$ we define
\begin{equation}\label{003}
w_{a,b}(x) = 
\bigl( \theta(a^2 x^{-2}) - 1 \bigr) \bigl( \theta(b^2 x^{-2}) -1 \bigr)
 - x^2 \bigl( \theta(a^2 x^2) - 1 \bigr) \bigl( \theta(b^2 x^2) -1 \bigr),
\end{equation}
where $\theta(x) = \sum_{ k \in {\Bbb Z}} e^{- \pi k^2 x}$ 
is the classical theta function. 
Using $w_{a,b}(x)$ the integrand $h_{\cal E}(x)$ in $\eqref{004}$ is expressed as 
\begin{equation}\label{005}
h_{\cal E}(x) = -{\frak e} \sum_{\nu} c(\nu) \int_{0}^{\infty} w_{a,\nu a^{-1}}(x) \frac{da}{a},  
\end{equation}
where $\frak e$ is a positive real constant and $c(\nu)$ are nonnegative real numbers given by  
\begin{equation} \label{006_1}
\sum_{\nu} c(\nu)\nu^{-s/2} = c_{\cal E}^{1-s} \zeta_{\cal E}(s)^2.
\end{equation}
See~\cite[sec.\,51]{Fe2} or~\cite[sec.\,8]{Fe3} for details. 
The asymptotic behavior of $h_{\cal E}(x)$ for small $x>0$ is given by  
\begin{equation}\label{006}
h_{\cal E}(e^{-t}) -  ( c_0 + c_1 t + c_2 t^2 + c_3 t^3 ) \to 0 \quad \text{as} \quad t \to +\infty
\end{equation}
for some constants $c_i$ $(0 \leq i \leq 3)$ with $c_3 \not=0$. 
Hence the remaining problem for the behavior of $h_{\cal E}(x)$ near $x=0$ 
is the behavior of the fourth derivative of $h_{\cal E}(e^{-t})$ for sufficiently large $t>0$.  

For $\nu >0$, we define
\begin{equation}\label{007}
V(x,\nu) = \int_{0}^{\infty} w_{a,\nu a^{-1}}(x) \frac{da}{a}
\end{equation}
and 
\begin{equation}\label{008}
Z(x,\nu) = \bigl( - x \frac{d}{dx} \bigr)^4 V(x,\nu). 
\end{equation}
Now we consider the series
\begin{equation}\label{009}
Z_{\cal E}(x) = \sum_{\nu} c(\nu) Z(x,\nu), 
\end{equation}
where $c(\nu)$ are the same in $\eqref{005}$. 
Since $ \frac{d}{dt} = -x \, \frac{d}{dx}$, 
we have $h_{\cal E}(e^{-t})^{\prime\prime\prime\prime} = - {\frak e}Z_{\cal E}(e^{-t})$. 
Thus the behavior of the fourth derivative of $h_{\cal E}(e^{-t})$ for large $t>0$ 
is obtained by the behavior of $Z_{\cal E}(x)$ for small $x>0$. 
Under the meromorphic continuation and the functional equation of $L(E,s)$, 
the relation between $Z_{\cal E}(x)$ and $L(E,s)$ are described as follows. 
\smallskip

\noindent
{\bf Theorem (Fesenko~\cite[Theorem 52]{Fe2}).}
{\it
Let $\cal E$ be a proper regular model of the elliptic curve $E/{\Q}$. 
Suppose that the model $\cal E$ is chosen as in section $42$ of ~\cite{Fe2}. 
Assume that}
\begin{enumerate}
\item[(F-1)] {\it $Z_{\cal E}(x)$ does not change its sign in some open interval $(0,x_0)$}, 
\item[(F-2)] {\it $L(E,s)$ has no real zeros in $(1,2)$}.  
\end{enumerate}
{\it Then all poles $\rho$ of $\zeta(s/2)\zeta(s)\zeta(s-1)/L(E,s)$ in the critical strip $0 < \Re(s) <2$ 
satisfy the Riemann hypothesis, namely, they lie on the line $\Re(\rho)=1$. 
}
\smallskip

This theorem is extended to more general situation (Theorem \ref{thm_05} in below). 
We often call (F-1) the {\it single sign property} of $Z_{\cal E}(x)$.
In a sense, Fesenko's theorem says that the single sign property of $Z_{\cal E}(x)$ 
is a sufficient condition for the Riemann hypothesis of $L(E,s)$. 
We show that it is a necessary condition under some technical assumption. 

\begin{theorem} \label{thm_001}
Suppose the Riemann hypothesis for $L(E,s)$. 
In addition, suppose that all zeros of $L({\cal E},s):=n_{\cal E}(s)^{-1}L(E,s)$ are simple 
except for the zero at $s=1$ and the estimate 
\begin{equation}
\sum_{0< \Im(\rho) \leq T}|L^\prime({\cal E},\rho)|^{-2} = O(T)
\end{equation}
holds, where $\rho$ runs all zeros of $L({\cal E},s)$ on the line $\Re(s)=1$. 
Then $Z_{\cal E}(x)$ is negative for sufficiently small $x>0$. 
\end{theorem}

This result is stated more precisely in section 4 (Theorem \ref{prop4}). 
\bigskip

We extend the above theorem of Fesenko in more general setting. 
Let ${\frak c}=\{c(\nu)\}$ be a sequence of nonnegative real numbers, 
and let $D_{\frak c}(s)$ be the Dirichlet series defined by 
\[
D_{\frak c}(s) = \sum_{\nu} c(\nu) \, \nu^{-s}.
\]
For the analytic treatment of $D_{\frak c}(s)$, 
we suppose the following three conditions for the nonnegative sequence ${\frak c} = \{c(\nu)\}$:
\begin{enumerate}
\item[(c-1)]  there exists a constant $M_{\varepsilon}>0$ such that 
$0 \leq c(\nu) \leq M_{\varepsilon} \, \nu^{\varepsilon}$ for any fixed $\varepsilon >0$, 
\item[(c-2)] there exists a constant $\eta_{\frak c}>0$ such that 
$D_{\frak c}(s)$ is continued holomorphically 
to the region 
\begin{equation}\label{300}
\sigma \geq 1 - \frac{\eta_{\frak c}}{\log(3+|t|)} \quad (s=\sigma+it)
\end{equation} 
except for the pole of order $\lambda_{\frak c} \geq 1$ at $s=1$, 
\item[(c-3)] there exist constants $M>0$ and $A>0$ such that
\begin{equation}
|D_{\frak c}(\sigma+it)| \leq M |t|^{A},
\end{equation}
whenever $\sigma \geq 1$ and $|t|>1$. 
\end{enumerate}
By (c-1) the abscissa of convergence $D_{\frak c}(s)$ is one. 
Thus $D_{\frak c}(s)$ has a singularity at $s=1$ by Landau's theorem. 
By (c-2) the singularity of $D_{\frak c}(s)$ at $s=1$ is constrained to be a pole. 
Condition (c-3) is a technical one, 
but it is ordinary satisfied by usual $L$-functions $L(s)$ and its reciprocal $1/L(s)$. 

Let $\gamma(s)$ be a meromorphic function on $\C$ satisfying
\begin{enumerate}
\item[($\gamma$-1)] $\gamma(s)$ is regular except for $s=1$, 
\item[($\gamma$-2)] $\gamma(s)$ has the pole of order $\lambda_\gamma \geq 1$ at $s=1$, 
\item[($\gamma$-3)] $\gamma(s)$ satisfies the uniform bound 
\[
|\gamma(\sigma+it)| \ll_{a,b,t_0} |t|^{-A} \quad (a \leq \sigma \leq b,~|t| \geq t_0)
\]  
for all real numbers $a \leq b$ and every real number $A>0$.
\end{enumerate}
We denote by $\kappa_\gamma(x)$ the inverse Mellin transform of $\gamma(s)$, namely, 
\[
\kappa_{\gamma}(x) = \frac{1}{2\pi i} \int_{(c)} \gamma(s) x^{-s} ds \quad (c>1).
\]
Then $\kappa_{\gamma}(x)$ is of rapid decay as $x\to+\infty$, 
namely, $\kappa_{\gamma}(x)=O(x^{-A})$ for all $A>0$ as $x\to+\infty$, 
and $\kappa_{\gamma}(x)=O(x^{-1-\delta})$ as $x \to 0^+$ for all $\delta>0$. 
\smallskip

We define
\begin{equation}\label{301}
h_{{\frak c},\gamma}^{\varepsilon,n}(x) = \sum_{\nu=1}^\infty c(\nu) V_{\gamma}^{\varepsilon,n}(x,\nu)
\end{equation}
for ${\frak c}=\{ c(\nu) \}$ satisfying (c-1), 
where $V_{\gamma}^{\varepsilon,n}(x,\nu)$ is the function defined by
\begin{equation} \label{301_1}
V_{\gamma}^{\varepsilon,n}(x,\nu) = \kappa_{\gamma}(\nu x^{-n}) - \varepsilon x^n \kappa_{\gamma}(\nu x^n) 
\end{equation}
for $\gamma(s)$ satisfying ($\gamma$-1), ($\gamma$-2), ($\gamma$-3) 
and $\varepsilon \in \{\pm 1\}$, $n, \nu \in {\Z}_{>0}$. 
The series on the right-hand side of \eqref{301} converges absolutely and uniformly 
on any compact set in $(0,\infty)$, since $\kappa_{\gamma}(x)$ is of rapid decay as $x \to +\infty$.     
The trivial functional equation 
$V_{\gamma}^{\varepsilon,n}(x^{-1},\nu) = - \varepsilon x^{-n} V_{\gamma}^{\varepsilon,n}(x,\nu)$ 
leads to the functional equation
\begin{equation}\label{302}
h_{{\frak c},\gamma}^{\varepsilon,n}(x^{-1}) = - \varepsilon x^{-n} h_{{\frak c},\gamma}^{\varepsilon,n}(x).
\end{equation}
Using $h_{{\frak c},\gamma}^{\varepsilon,n}(x)$ we define
\begin{equation}\label{303}
\omega_{{\frak c},\gamma}^{\varepsilon,n}(s)=\int_{0}^{1} h_{{\frak c},\gamma}^{\varepsilon,n}(x) \, x^{s-n} \frac{dx}{x}. 
\end{equation} 
and 
\begin{equation}\label{304}
Z_{{\frak c},\gamma}^{\varepsilon,n}(x) 
= \Bigl( -x\frac{d}{dx} \Bigr)^{\lambda_{\frak c}+\lambda_\gamma} h_{{\frak c},\gamma}^{\varepsilon,n}(x), 
\end{equation}
where $\lambda_{\frak c}$ is the order of the pole of $D_{\frak c}(s)$ at $s=1$ 
and $\lambda_{\gamma}$ is the order of the pole of $\gamma(s)$ at $s=1$.  
These functions $h_{{\frak c},\gamma}^{\varepsilon,n}(x)$, $\omega_{{\frak c},\gamma}^{\varepsilon,n}(s)$ 
and $Z_{{\frak c},\gamma}^{\varepsilon,n}(x)$ 
are a generalization of $h_{\cal E}(x)$, $\omega_{\cal E}(s)$ and $Z_{\cal E}(x)$, respectively. 
In fact, if we take $D_{\frak c}(s)=c_{\cal E}^{1-s}\zeta_{\cal E}(2s)^2$, 
$\gamma(s)=\widehat{\zeta}(s)^2$, $\varepsilon=+1$ and $n=2$, then we get them. 
In general, when we study $h_{\cal E}(x)$ of the model $\cal E$ of the elliptic curve over the algebraic number field $k$, 
we need the case $\gamma(s)=\prod_{1\leq i \leq I}\widehat{\zeta}_{k_i}(s)^2$, 
where $k_i$ are finite extensions of $k$ which include $k$ itself.  
\smallskip

The single sign property of 
$Z_{{\frak c},\gamma}^{\varepsilon, n}(s)$ 
implies the nonexistence of poles of $\omega_{{\frak c},\gamma}^{\varepsilon, n}(s)$ 
around the line $\Re(s)=n$ except for $s=n$. 
\begin{theorem}\label{thm_05}
Let ${\frak c}$ be a nonnegative sequence satisfying (c-1), (c-2) and (c-3). 
Suppose that there exists $x_0>0$ such that $Z_{{\frak c},\gamma}^{\varepsilon, n}(s)$ has a single sign on $(0,x_0)$. 
Then there exists $\delta>0$ such that 
$\omega_{{\frak c},\gamma}^{\varepsilon, n}(s)$  is continued holomorphically to the right-half plane 
$\Re(s)> n -\delta$ except for the pole $s=n$. 

Further suppose that $\omega_{{\frak c},\gamma}^{\varepsilon, n}(s)$ 
is continued meromorphically to the right-half plane $\Re(s)>\sigma_0$ 
for some $\sigma_0 < n$ without poles on the open interval $(\sigma_0, n)$. 
Then $\omega_{{\frak c},\gamma}^{\varepsilon, n}(s)$ has no pole in the vertical strip $ \sigma_0 < \Re(s) < n$. 
\end{theorem}
The case $\lambda_{\frak c}=2$ and $\gamma(s) = \widehat{\zeta}(s)^2~(\lambda_\gamma=2)$ 
is essentially Fesenko's result in the above. 
Theorem \ref{thm_05} is proved in section 2.  
Now we consider the problem:
\smallskip 

\noindent
{\bf Problem.} {\it 
For which kind of ${\frak c}$, 
will $Z_{{\frak c},\gamma}^{\varepsilon,n}(x)$ keep its sign for sufficiently small $x>0$ ?}
\smallskip

\noindent
This is a question on the property ($\ast$) in ~\cite[51 in section 4.3]{Fe2}. 
Unfortunately, theoretical progress on the problem is not yet obtained. 
However, we would give some remark on the problem in the final section.  
\smallskip

Now we back to the case of $\zeta_{\cal E}(s)$. 
For a concrete elliptic curve $E/{\Q}$ with a small conductor, 
we can see if (F-2) is supported by computations . 
In fact (F-2) holds whenever the conductor of $E/{\Q}$ is less than $8000$ (see Rubinstein~\cite{Ru}). 
Hence, for such $E/{\Q}$, our interest is in the behavior of $Z_{\cal E}(x)$ for small $x>0$. 
Recall the formula \eqref{009} of $Z_{\cal E}(x)$. 
As a first step of the research for $Z_{\cal E}(x)$, 
we give a series expansion of $Z(x,\nu)$ consisting of $K$-Bessel functions: 
\[
Z_\nu(x) = 4 \Bigl( -x\frac{d}{dx} \Bigr)^4 \Bigl(
 \sum_{N=1}^\infty \sigma_0(N) \bigl(  K_0(2 \pi N \nu x^{-2}) - x^2 K_0(2 \pi N \nu x^{2}) \bigr)
\Bigr),
\]
where $\sigma_0(N)=\sum_{d|N}1$ and $K_0(t)$ is the $K$-Bessel function of index $0$. 
Using this expansion, we find that each $Z(x,\nu)$ is negative for sufficiently small $x>0$. 
Also, it enables us to see if (F-1) is supported by computations for given concrete elliptic curve $E/{\Q}$ 
whenever the conductor is small (see the table \cite{Fe4}). 
These results and other results relating $Z_{\cal E}(s)$ are stated and proved in section 3. 
Of course, they do not ensure (F-1). 
Further study of (F-1) will be conducted in the future. 
One attractive approach for (F-1) is 
to study $h_{\cal E}(x)$ using properly the detail structure of the boundary of a two-dimensional adelic object 
appearing in the integral representation of $h_{\cal E}(x)$ and $\omega_{\cal E}(s)$ 
(see \cite[section 4.3]{Fe2} and \cite[section 6]{Fe3}).  
\smallskip

Finally, in section 5, we remark on the single sign property of $Z(x)$ 
from a viewpoint of an Euler product. 
In the section, we study $Z(x)$ corresponding to the Dirichlet series
\[
D(s) = \frac{\zeta(2s)^2 \zeta(2s-1)^2}{L(2s-1/2)^2}, 
\]
where $L(s)$ is a function defined by Euler products of degree two. 
We observe that we can obtain the best possible estimate of $Z(x)$ for small $x>0$ for ``almost all $L(s)$'',  
but that $Z(x)$ may have oscillation near $x=0$ in general. 
Hence we set our interest on some special class of $L$-functions. 
We choose the Selberg class as such class of $L$-functions, 
and state a conjecture for the single sign property of $Z(x)$ 
corresponding to the $L$-functions in the Selberg class. 
Interestingly, the study of such $Z(x)$ is related not only to 
the Riemann hypothesis of $L(E,s)$ but also to the BSD-conjecture. 
In fact two dimensional adelic analysis programme 
suggests a new method to prove the rank part of the BSD (\cite[section 9]{Fe3}). 
In section 5.3, we study $Z(x)$ from the viewpoint of the ``partial Euler product'' of $L(E,s)$ 
according to Goldfeld~\cite{MR679556} and Conrad~~\cite{MR2124918}. 

\smallskip

\noindent
{\bf Notations.}
We always denote by $\varepsilon$ an arbitrary small positive real number. 
For a positive valued function $g(x)$, 
we use Landau's $f(x)=O(g(x))$ and Vinogradov's $f(x) \ll g(x)$ as the same meaning. 
More precisely $f(x)=O(g(x))$ or $f(x) \ll g(x)$ for $x \in X$ means that 
$|f(x)| \leq C g(x)$ for any $x \in X$ and some constant $C \geq 0$. 
Any value $C$ for which this holds is called an implied constant. 
Since a constant is often a function depending on a variable, 
the ``implied constant'' will sometimes depends on other parameters, 
which we explicitly mention at important points. 
Also we use Landau's $f(x)=o(g(x))$ for $x \to x_0$ 
in the meaning that 
for any $\varepsilon>0$ there exists a neighborhood $U_{\varepsilon}$ of $x_0$ 
such that $|f(x)| \leq \varepsilon g(x)$ for any $x \in U_{\varepsilon}$. 
\smallskip

\noindent{\bf Acknowledgment.} 
The author thanks Ivan Fesenko for his much encouragement and many helpful comments to this research. 
The author thanks the first referee of 
the paper for his suggestions to 
improve and simplify the proof of several results. 

\section{Boundary term for general Dirichlet series}

In this part, we describe the relation between the Dirichlet series 
$D_{\frak c}(s)$ and the boundary term $\omega_{{\frak c},\gamma}^{\varepsilon,n}(s)$. 
Throughout this section, we assume that the nonnegative sequence 
${\frak c}$ satisfies (c-1), (c-2) and (c-3) and $\gamma(s)$ satisfies ($\gamma$-1), ($\gamma$-2) and ($\gamma$-3) 
whenever we do not mention a condition for $\frak c$ or $\gamma(s)$. 


\begin{proposition} \label{prop_01}
Let $\frak c$ be a nonnegative sequence satisfying (c-1). 
Then 
\begin{equation}\label{0201}
n^{-1} \gamma(s/n) D_{\frak c}(s/n) = \xi(s) + \varepsilon \, \xi(n-s) - \varepsilon \, \omega_{{\frak c},\gamma}^{\varepsilon,n}(s)
\end{equation}
for $\Re(s)>n$, where $\omega_{{\frak c},\gamma}^{\varepsilon,n}(s)$ is defined by \eqref{303} 
and $\xi(s)$ is an entire function given by 
\[
\xi(s)=\int_{1}^{\infty} \Bigl(~ \sum_{\nu=1}^\infty c(\nu) \kappa_{\gamma}(\nu x^{n})  ~\Bigr) x^{s} \frac{dx}{x}
\]
using the inverse Mellin transform $\kappa_\gamma(x)$ of $\gamma(s)$.
\end{proposition}
This is a generalization of the well known integral representation 
\[
\widehat{\zeta}(s) = 
\int_{1}^{\infty} (\theta(x^2)-1) \, x^{s} \frac{dx}{x} 
+ \int_{1}^{\infty} (\theta(x^2)-1) \, x^{1-s} \frac{dx}{x} 
+ \Bigl( \frac{1}{s-1} - \frac{1}{s} \Bigr).
\]
In fact, if we take ${\frak c}=\{c(\nu)=1\}$, $\gamma(s)=\pi^{-s/2}\Gamma(s/2)$, $\varepsilon=+1$ and $n=1$, 
then $\sum_{\nu=1}^\infty c(\nu)\kappa_{\gamma}(\nu x)=\theta(x^2)-1$ 
and $\omega_{{\frak c},\gamma}^{+1,1}(s)=(1-s)^{-1}+s^{-1}$. 
If we take $D_{\frak c}(s)=c_{\cal E}^{1-s}\zeta_{\cal E}(2s)^2$, 
$\gamma(s)=\widehat{\zeta}(s)^2$, $\varepsilon=+1$ and $n=2$, 
we obtain $\omega_{\cal E}(s)$. In general, 
to study the boundary term $\omega_{\cal E}(s)$ over the algebraic number field $k$, 
we need the case $\gamma(s)=\prod_{1\leq i \leq I}\widehat{\zeta}_{k_i}(s)^2$, 
where $k_i$ are finite extension of $k$ which include $k$ itself.  

\begin{corollary} 
Let ${\frak c}=\{c(\nu)\}$ be a nonnegative sequence satisfying (c-1), and let $\delta>0$.   
Then the meromorphic continuation of $D_{\frak c}(s)$ to the right-half plane $\Re(s) > 1 - \delta$ 
is equivalent to the meromorphic continuation of $\omega_{{\frak c},\gamma}^{\varepsilon,n}(s)$ 
to  the right-half plane $\Re(s)>n(1-\delta)$. 
In addition, the functional equations $\gamma(s)D_{\frak c}(s) = \varepsilon \, \gamma(1-s)D_{\frak c}(1-s)$ 
and $\omega_{{\frak c},\gamma}^{\varepsilon,n}(s)= \varepsilon \, \omega_{{\frak c},\gamma}^{\varepsilon,n}(n-s)$ 
are equivalent to each other. 
\end{corollary}
The form of $\omega_{{\frak c},\gamma}^{\varepsilon,n}(s)$ is completely determined, 
if we assume the functional equation of the product $\gamma(s)D_{\frak c}(s)$. 
\begin{proposition} \label{prop_02}
Suppose that $\gamma(s)D_{\frak c}(s)$ is continued meromorphically to $\C$ 
and satisfies the functional equation 
$\gamma(s)D_{\frak c}(s) = \varepsilon \gamma(1-s) D_{\frak c}(1-s)$ 
for $\varepsilon \in \{\pm 1\}$. 
Further suppose that there exists a real number $A \geq 0$ 
and a strictly increasing sequence of positive real numbers $\{t_n\}_{n \geq 1}$ 
such that 
$|D_{\frak c}(\sigma+it_n)| \ll t_n^A$ uniformly for $\sigma \in [-1/2,5/2]$. 
Then we have
\begin{equation}\label{p2}
-\omega_{{\frak c},\gamma}^{\varepsilon,n}(s) = 
  \sum_{m=1}^{\lambda_{\frak c}+\lambda_\gamma} C_m \Bigl( \frac{1}{(s-n)^m}  + \varepsilon \frac{(-1)^m}{s^m} \Bigr) 
+ \sum_{\rho}\sum_{m=1}^{m_\rho} \frac{C_{\rho,m}}{(s - \rho)^m},
\end{equation}  
where constants $C_m$ are given by 
\begin{equation}
n^{-1} \gamma(s/n) D_{\frak c}(s/n) = \sum_{m=1}^{\lambda_{\frak c}+\lambda_\gamma} \frac{C_m}{(s-n)^m} + O(1)
\quad \text{as} \quad s \to n,
\end{equation}
the sum $\sum_\rho$ runs over all poles $\rho$ of $\gamma(s/n) D_{\frak c}(s/n)$ 
satisfying $0 < \Re(\rho) < n$, and the constant $C_{\rho,m}$ is given by
\begin{equation}
n^{-1} \gamma(s/n) D_{\frak c}(s/n) = \sum_{m=1}^{m_\rho} \frac{C_{\rho,m}}{(s-\rho)^m} + O(1)
\quad \text{as} \quad s \to \rho.
\end{equation}
Moreover, the sum $\sum_{\rho}$ in \eqref{p2} 
converges uniformly on every compact subset in $\C$. 
\end{proposition}
We define
\begin{equation} \label{0205}
Z_{\gamma}^{\varepsilon,n}(x,\nu)=
\Bigl(-x\frac{d}{dx} \Bigr)^{\lambda_{\frak c}+\lambda_{\gamma}} V_{\gamma}^{\varepsilon,n}(x,\nu). 
\end{equation}
Note that
\[
Z_{{\frak c},\gamma}^{\varepsilon,n}(x,\nu) = \sum_{\nu=1}^\infty c(\nu) Z_{\gamma}^{\varepsilon,n}(x,\nu). 
\]
by \eqref{304}. 
In general, $Z_{\gamma}^{\varepsilon,n}(x,\nu)$ has the following dilation formula. 
\begin{proposition}~\label{thm4} We have
\begin{equation}\label{107}
Z_{\gamma}^{\varepsilon,n}(x,\nu) 
= \frac{\nu_0}{\nu} \, Z_{\gamma}^{\varepsilon,n}\bigl(x (\nu/\nu_0)^{1/n}, \, \nu_0 \,\bigr) 
+ R(x,\nu ; \nu_0).  
\end{equation}
For any fixed $\delta>0$, the second term $R(x, \nu ; \nu_0)$ satisfies the estimate 
\[
R(x, \nu ; \nu_0) \ll \nu^{-1-\delta}x^{n(1+\delta)}
 + \nu^{-1} \nu_0^{-\delta}\bigl( x (\nu/\nu_0)^{1/n} \bigr)^{n(1+\delta)}
\]
for $0<x<1$, where the implied constant depends only on $\gamma$, $n$ and $\delta$. 
In particular, 
\[
R(x, \nu ; \nu_0 ) \ll \nu^{-1} \, \nu_0^{-\delta} R^{n(1+\delta)}
\]
under the condition $\nu > \nu_0$ and  $x (\nu/\nu_0)^{1/n} \leq R <1$ for given $R>0$. 
\end{proposition}
Further we define
\begin{equation} \label{0207}
\aligned
Z_{\gamma,0}^{n}(x,\nu) & = 
\Bigl(-x\frac{d}{dx} \Bigr)^{\lambda_{\frak c}+\lambda_{\gamma}} \bigl( x^n \kappa_{\gamma}(\nu x^n) \bigr), \\
Z_{\gamma,1}^{n}(x,\nu) & = 
\Bigl(-x\frac{d}{dx} \Bigr)^{\lambda_{\frak c}+\lambda_{\gamma}} \bigl( \kappa_{\gamma}(\nu x^{-n}) \bigr)
\endaligned
\end{equation}
and 
\begin{equation} \label{0208}
Z_{{\frak c},\gamma,i}^{n}(x) = \sum_{\nu=1}^\infty c(\nu) Z_{\gamma,i}^{n}(x,\nu). 
\quad (i=0,1)
\end{equation}
By \eqref{301_1}, \eqref{302} and \eqref{304}, we have
\[
Z_{{\frak c},\gamma}^{\varepsilon, n}(x,\nu) = Z_{{\frak c},\gamma,1}^{n}(x,\nu) - \varepsilon Z_{{\frak c},\gamma,0}^{n}(x,\nu).
\] 
\begin{proposition} \label{prop_004}
We have
\[
Z_{{\frak c},\gamma}^{\varepsilon, n}(x,\nu) = 
- \varepsilon Z_{{\frak c},\gamma,0}^{n}(x,\nu) 
+ O(x^A)
\]
as $x \to 0^+$ for any fixed real number $A$, 
where the implied constant depends only on ${\frak c}$, $\gamma(s)$, $n$ and $A$
\end{proposition}
\subsection{Proof of Proposition \ref{prop_01}} 
We put
\begin{equation} \label{305}
f_{{\frak c},\gamma}^{n}(x) = \frac{1}{2\pi i}\int_{\Re(s)=c} \gamma(s) D_{\frak c}(s) x^{ns} ds
\end{equation}
for $c>1$. 
The conditions for $\frak c$ and $\gamma(s)$ imply that  $f_{{\frak c},\gamma}^{n}(x)$ is well-defined 
and the Mellin inversion formula gives 
\begin{equation} \label{306}
n^{-1} \gamma(s/n) D_{\frak c}(s/n) = \int_{0}^{\infty} f_{{\frak c},\gamma}^{n}(x^{-1}) x^{s} \frac{dx}{x}
\end{equation}
for $\Re(s)>n$. On the right-hand side, we have 
\[
\aligned
\int_{0}^{\infty} f_{{\frak c},\gamma}^{n}(x^{-1}) x^{s} \frac{dx}{x}
& = \int_{1}^{\infty} f_{{\frak c},\gamma}^{n}(x^{-1}) x^{s} \frac{dx}{x} 
+ \varepsilon \int_{1}^{\infty} f_{{\frak c},\gamma}^{n}(x^{-1}) x^{n-s} \frac{dx}{x} \\
& \quad - \varepsilon \int_{0}^{1} \bigl( 
f_{{\frak c},\gamma}^{n}(x) - \varepsilon x^n  f_{{\frak c},\gamma}^{n}(x^{-1}) 
\bigr) x^{s-n} \frac{dx}{x}
\endaligned
\]
for $\varepsilon \in \{\pm 1\}$. Using the Dirichlet series expansion of $D_{\frak c}(s)$, we have
\begin{equation*}
f_{{\frak c},\gamma}^{n}(x) = \sum_{\nu=1}^\infty c(\nu) \kappa_\gamma(\nu x^{-n}). 
\end{equation*}
Therefore
\[
f_{{\frak c},\gamma}^{n}(x) - \varepsilon x^n  f_{{\frak c},\gamma}^{n}(x^{-1}) 
= \sum_{\nu=1}^\infty c(\nu) \bigl( \kappa_\gamma(\nu x^{-n}) - \varepsilon x^n \kappa_\gamma(\nu x^n) \bigr) 
= h_{{\frak c},\gamma}^{\varepsilon,n}(x)
\]
by definition \eqref{301} and \eqref{301_1}, 
and 
\begin{equation} \label{307}
\aligned
\int_{0}^{\infty} f_{{\frak c},\gamma}^{n} & (x^{-1}) x^{s} \frac{dx}{x} \\
& = \int_{1}^{\infty} f_{{\frak c},\gamma}^{n}(x^{-1}) x^{s} \frac{dx}{x} 
+ \varepsilon \int_{1}^{\infty} f_{{\frak c},\gamma}^{n}(x^{-1}) x^{n-s} \frac{dx}{x} 
- \varepsilon \, \omega_{{\frak c},\gamma}^{\varepsilon,n}(s)
\endaligned
\end{equation}
by definition \eqref{303}. 
By \eqref{306} and \eqref{307}, we obtain \eqref{0201} for $\Re(s)>n$. 
\hfill $\Box$ 

\subsection{Proof of Proposition \ref{prop_02}}
By definition $\eqref{305}$ of $f_{{\frak c},\gamma}^n(x)$,  we have
\[
f_{{\frak c},\gamma}^n(x) = \frac{1}{2 \pi i} \int_{Re(s)=c}
n^{-1} \gamma(s/n) D_{\frak c}(s/n) \, x^{s} ds \quad (c>n).
\]
Moving the path of integration to the line $\Re(s)=c^\prime$ with $c^\prime<0$, we have 
\[
\aligned
f_{{\frak c},\gamma}^n(x)
& =  
 \sum_{m=1}^{\lambda_{\frak c}+\lambda_\gamma} \frac{C_{m}}{(m-1)!} (x^n + (- 1)^m) (\log x)^{m-1} 
+ \sum_{\rho} x^{\rho} \sum_{m=1}^{m_\rho} \frac{C_{\rho,m}}{(m-1)!} (\log x)^{m-1} \\
& \quad + \frac{1}{2 \pi i} \int_{\Re(s)=c^\prime}
n^{-1} \gamma(s/n) D_{\frak c}(s/n) \, x^{s} ds.
\endaligned
\]
This process is justified by ($\gamma$-3) 
and the assumption for $D_{\frak c}(s)$ in the proposition. 
Using the functional equation $\gamma(s)D_{\frak c}(s)= \varepsilon \gamma(1-s)D_{\frak c}(1-s)$, 
the third integral of the right-hand side equals 
\[
\frac{1}{2 \pi i} \int_{\Re(s)=n-c^\prime}
\varepsilon \, n^{-1} \gamma(s/n) D_{\frak c}(s/n) \, x^{n-s} ds
= \varepsilon x^n  f_{{\frak c},\gamma}^{n}(x^{-1}).
\]
Hence we have
\[
\aligned
h_{{\frak c},\gamma}^{\varepsilon,n}(x) 
& = f_{{\frak c},\gamma}^{n}(x) - \varepsilon x^n  f_{{\frak c},\gamma}^{n}(x^{-1})  \\ 
& = 
\sum_{m=1}^{\lambda_{\frak c}+\lambda_\gamma} \frac{C_{m}}{(m-1)!} (x^n + (-1)^m) (\log x)^{m-1} 
+ \sum_{\rho} x^{\rho} \sum_{m=1}^{m_\rho} \frac{C_{\rho,m}}{(m-1)!} (\log x)^{m-1}.
\endaligned
\]
This implies $\eqref{p2}$ by definition \eqref{303} of $\omega_{{\frak c},\gamma}^{\varepsilon,n}(s)$, 
since if $\rho$ is a pole of $D_{\frak c}(s/n)$ 
in $0 < \Re(s) <n$ then $n-\rho$ is also a pole of $D_{\frak c}(s/n)$ 
which has the same multiplicity of $\rho$ and $C_{n-\rho,m}(s) = (-1)^m \, C_{\rho,m}$.
\hfill $\Box$

\subsection{Proof of Theorem \ref{thm_05}}
To prove Theorem \ref{thm_05}, we prepare the lemma 
for the asymptotic behavior of $f_{{\frak c},\gamma}^n(x)$ 
and an analogue of Landau's theorem for Laplace transforms.
\begin{lemma} \label{lem_01}
Suppose that
\begin{equation}\label{312_1}
n^{-1}\gamma(s/n)D_{\frak c}(s/n) 
= \frac{C_{\lambda_{\frak c}+\lambda_{\gamma}}}{(s-n)^{\lambda_{\frak c}+\lambda_{\gamma}}} + \dots  + \frac{C_1}{s-n} + O(1) 
\quad (C_{\lambda_{\frak c}+\lambda_{\gamma}} \not = 0)
\end{equation}
in a neighborhood of $s = n$. Then 
\begin{equation}\label{312_2}
f_{{\frak c},\gamma}^n(x) = x^n \sum_{m=0}^{\lambda_{\frak c}+\lambda_{\gamma}-1} \frac{C_{m+1}}{m!} \, (\log x)^m + o(x^n). 
\end{equation}
for $x>1$. 
\end{lemma}
\begin{lemma}
Let $f(t)$ be a real valued function on $(0,\infty)$. 
Suppose that there exists $t_0>0$ such that $f(t)$ does not change its sign for $t>t_0$ 
and the abscissa $\sigma_c$ of the convergence of the integral 
\[
F(s)=\int_{0}^{\infty} f(t) e^{-st} dt
\]
is finite. 
Then $F(s)$ has a singularity on the real axis at the point $s=\sigma_c$.
\end{lemma}
\begin{pf}
Refer to section $5$ of chapter II in~\cite{MR0005923}.  
\end{pf}

\noindent
{\bf Proof of Lemma \ref{lem_01}.} 
We denote $n^{-1} \gamma(s/n) D_{\frak c}(s/n)$ by $F(s)$, 
and $f_{{\frak c},\gamma}^n(x)$ by $f(x)$. 
Note that $f(x)=\frac{1}{2\pi i}\int_{\Re(s)=c} F(x) x^{s} ds~(c>n)$. 
By condition (c-2), $F(s)$ is extended holomorphically to the region $\eqref{300}$ 
except for the pole of order $\lambda_{\frak c} + \lambda_\gamma$ at $s=n$. 
We put
\[
P(s) = \frac{C_{\lambda_{\frak c} + \lambda_\gamma}}{(s-n)^{\lambda_{\frak c} + \lambda_\gamma}} + \dots 
+ \frac{C_2}{(s-n)^2} + C_1 \bigl(\frac{1}{s-n} - \frac{1}{s-n+1} \bigr)
\]
and 
\[
p(x) = x^n \sum_{m=1}^{\lambda_{\frak c} + \lambda_\gamma-1} \frac{C_{m+1}}{m!}(\log x)^m +C_1 x^{n-1}(x-1). 
\]
Then $F(s)-P(s)$ is holomorphic for $\Re(s) \geq n$. 
Since 
\[
p(x)=x^n \sum_{m=1}^{\lambda_{\frak c} + \lambda_\gamma-1} \frac{C_{m+1}}{m!}(\log x)^m +o(x^n),
\]
it suffices to show that $(f(x) - p(x))/x^n = o(1)$ as $x \to \infty$. 
Using
\[
\frac{x^n (\log x)^m}{m!} =
\frac{1}{2\pi i} \int_{c - i \infty}^{c + i \infty} \frac{x^{s}}{(s-n)^{m+1}} ds \quad (x>1,c>n)
\]
for $m=0,1,2,\dots$ and
\[
x^{n-1}(x-1) = \frac{1}{2\pi i} \int_{c - i \infty}^{c + i \infty} 
\frac{x^{s}}{(s-n+1)(s-n)} ds \quad (x>1, c>n), 
\]
we have
\[
f(x) - p(x) = 
\frac{1}{2 \pi i} \int_{c - i \infty}^{c + i \infty} 
( F(s) - P(s) ) x^{s} ds \quad (x>1,c>n).  
\]
By (c-3), we can move the path of integration to the line $\Re(s)=n$, and get
\[
\frac{f(x) - p(x)}{x^n} = \frac{1}{2 \pi i} 
\int_{-\infty}^{\infty} ( F(n+it) - P(n+it) ) e^{it \log x} dt \quad (x>1).
\]
Since $F(s) - P(s)$ is holomorphic on $\Re(s)=n$, 
(c-3) and ($\gamma$-3) give  
\[
\int_{-\infty}^{\infty} | F(n+it) - P(n+it) | dt < \infty. 
\]
The Riemann-Lebesgue lemma in the theory 
of Fourier series states that 
\[
\lim_{x \to + \infty} \int_{-\infty}^{\infty} f(t) e^{itx} dt =0
\]
if the integral $\int_{-\infty}^{\infty} |f(t)| dt$ converges. 
Thus we obtain 
\[
\lim_{x \to +\infty} \int_{-\infty}^{\infty} ( F(n+it) - P(n+it) ) e^{it \log x} dt=0.
\]
This implies $(f(x) - p(x))/x^n=o(1)$ as $x \to +\infty$. 
\hfill $\Box$
\bigskip

\noindent
{\bf Proof of Theorem \ref{thm_05}.}
We denote $h_{{\frak c},\gamma}^{\varepsilon,n}(e^{-t})$ by $H(t)$. 
Then our single sign assumption for $Z_{{\frak c},\gamma}^{\varepsilon,n}(x)$ implies 
$H^{(\lambda_{\frak c}+\lambda_{\gamma})}(t)$ does not change its sign for $t>t_0$.  
We have
\[
H(t)
= h_{{\frak c},\gamma}^{\varepsilon,n}(e^{-t}) 
= - \sum_{m=0}^{\lambda_{\frak c}+\lambda_{\gamma}-1} \frac{C_{m+1}}{m!} t^{m} + o(1)
\]
as $t \to +\infty$ by $h_{{\frak c},\gamma}^{\varepsilon,n}(x)= f_{{\frak c},\gamma}^n(x) - \varepsilon x^n f_{{\frak c},\gamma}^n(x^{-1})$ 
and Lemma \ref{lem_01}, where $C_m$ are numbers in $\eqref{312_1}$. 
Hence we have
\[
\omega_{{\frak c},\gamma}^{\varepsilon,n}(s+n)
 = \int_{0}^{\infty} H(t) e^{-st} dt = 
-\sum_{m=0}^{\lambda_{\frak c}+\lambda_{\gamma}-1}\frac{H^{(m)}(0)}{s^{m+1}} 
+ \frac{1}{s^{\lambda_{\frak c}+\lambda_{\gamma}}} \int_{0}^{\infty} H^{(\lambda_{\frak c}+\lambda_{\gamma})}(t) e^{-st} dt.
\]
Thus we obtain
\begin{equation}\label{313}
\int_{0}^{\infty} H^{(\lambda_{\frak c}+\lambda_{\gamma})}(t) e^{-st} dt
= s^{\lambda_{\frak c}+\lambda_{\gamma}} \Bigl(\omega_{{\frak c},\gamma}^{\varepsilon,n}(s+n) + \sum_{m=1}^{\lambda_{\frak c}+\lambda_{\gamma}} \frac{C_m}{s^m} \Bigr).
\end{equation}
Here definition $\eqref{312_1}$ of $C_m$ implies that the right-hand side is regular around $s=0$, 
since $\omega_{{\frak c},\gamma}^{\varepsilon,n}(s)=\xi(s)+\xi(2-s)- n^{-1}\gamma(s/n)D_{\frak c}(s/n)$ 
with the entire function $\xi(s)$.  

Let $\sigma_c$ be the abscissa of convergence of the integral 
\[
{\cal H}(s)=\int_{0}^{\infty} H^{(\lambda_{\frak c}+\lambda_{\gamma})}(t) e^{-st} dt.
\]
Then ${\cal H}(s)$ has a singularity on the real axis at $s=\sigma_c$, 
since $H^{(\lambda_{\frak c}+\lambda_{\gamma})}(t)$ is real-valued and does not change its sign for $t>t_0$. 
Thus the abscissa $\sigma_c$ must be negative, namely, 
$\sigma_c=-\delta$ for some $\delta>0$. 
In particular ${\cal H}(s-n)$ is regular for $\Re(s)>n-\delta$. 
Hence the function $\omega_{{\frak c},\gamma}^{\varepsilon,n}(s)$ is continued holomorphically 
to the right-half plane $\Re(s)>n-\delta$ except for $s=n$,  
because of  
\begin{equation}\label{315}
\omega_{{\frak c},\gamma}^{\varepsilon,n}(s)
=-\sum_{m=1}^{\lambda_{\frak c}+\lambda_{\gamma}} \frac{C_m}{(s-n)^m} + \frac{{\cal H}(s-n)}{(s-n)^{\lambda_{\frak c}+\lambda_{\gamma}}}.
\end{equation}
Now we suppose that $\omega_{{\frak c},\gamma}^{\varepsilon,n}(s)$ is continued meromorphically to $\Re(s)>\sigma_0$ 
with no real pole on $(\sigma_0,n)$. 
Then, by $\eqref{313}$, ${\cal H}(s)$ is continued to $\Re(s)>\sigma_0-n$ and has no pole on $(\sigma_0-n,\infty)$. 
This implies $\sigma_c \leq \sigma_0-n$ and that ${\cal H}(s-n)$ is regular for $\Re(s)>\sigma_0$. 
Hence, by $\eqref{315}$, $\omega_{{\frak c},\gamma}^{\varepsilon,n}(s)$  has no pole for $\Re(s) > \sigma_0$ except for $s=n$. 
\hfill $\Box$

\subsection{Proof of Proposition \ref{thm4}}
By definition \eqref{301_1}, we have
\begin{equation} \label{0320}
V_{\gamma}^{\varepsilon,n}(x,\nu) = \frac{1}{2 \pi i}\int_{(c)} 
 \nu^{-s} (x^{ns} - \varepsilon x^{n(1-s)}) \, \gamma(s) ds 
\quad (c>1).
\end{equation}
For $\nu_0>0$ we define  
\begin{equation*}
R(x,\nu ; \nu_0) = Z_{\gamma}^{\varepsilon,n}(x,\nu) 
- \frac{\nu_0}{\nu} \, Z_{\gamma}^{\varepsilon,n}\bigl(x (\nu/\nu_0)^{1/n}, \, \nu_0 \,\bigr) .
\end{equation*}
By definition \eqref{0205} of $Z_{\gamma}^{\varepsilon,n}(x,\nu)$, we have 
\begin{equation*}
Z_{\gamma}^{\varepsilon,n}(x,\nu) = (-1)^{\lambda_{\frak c}+\lambda_{\gamma}} \frac{1}{2 \pi i}\int_{(c)} 
 \nu^{-s} ((ns)^{\lambda_{\frak c}+\lambda_{\gamma}} x^{ns} 
- \varepsilon (n(1-s))^{\lambda_{\frak c}+\lambda_{\gamma}} x^{n(1-s)}) \, \gamma(s) ds 
\end{equation*}
for $c>1$. 
Hence 
\begin{equation*}
R  (x,\nu ; \nu_0) = (-1)^{\lambda_{\frak c}+\lambda_{\gamma}}
\frac{1}{2 \pi i}\int_{(c)} 
x^{ns} ( \nu^{-s} -  \nu^{s-1}\nu_0^{1-2s}) \, (ns)^{\lambda_{\frak c}+\lambda_{\gamma}} \gamma(s) ds
\quad (c>1).
\end{equation*}
Moving the path of integration to the line $\Re(s)=1+\delta~(\delta>0)$, we obtain
\[
R  (x,\nu ; \nu_0) \ll \nu^{-1-\delta}x^{n(1+\delta)} + \nu^{\delta}\nu_0^{-1-2\delta}x^{n(1+\delta)}
\]
where the implied constant depends only on $\delta>0$. 
This shows the first estimate in Proposition \ref{thm4}, 
since $\nu^{\delta}\nu_0^{-1-2\delta}x^{n(1+\delta)} 
= \nu^{-1}\nu_{0}^{-\delta}[x (\nu/\nu_0)^{1/n}]^{n(1+\delta)}$. 
\hfill $\Box$

\subsection{Proof of Proposition \ref{prop_004}} 
It suffices to show that $Z_{\gamma,1}^{n}(x) = O(x^A)$ as $x \to 0^+$. 
By definition \eqref{301_1} and \eqref{0207}, we have
\[
\aligned
Z_{\gamma,1}^{n}(x) 
& = \Bigl( -x\frac{d}{dx} \Bigr)^{\lambda_{\frak c}+\lambda_\gamma} 
\frac{1}{2\pi i} \int_{(c)} \gamma(s)  D_{\frak c}(s) x^{ns} ds \\
& = \frac{1}{2\pi i} \int_{(c)}  \gamma(s)  D_{\frak c}(s) (-ns)^{\lambda_{\frak c}+\lambda_\gamma} x^{ns} ds 
\quad (c>1).
\endaligned
\]
Moving the path of integration to the right, we obtain the above assertion. 
\hfill $\Box$


\section{Case of two-dimensional zeta integral over $\Q$}

In this part, we apply the results in section 2 to $h_{\cal E}(x)$, $Z_{\cal E}(x)$ and $\omega_{\cal E}(s)$, 
and add more detailed consideration.  
More precisely, we study the function $V_{\gamma}^{\varepsilon,n}(x,\nu)$ in \eqref{301_1}, 
its derivative and $Z_{{\frak c},\gamma}^{\varepsilon,n}(x)$ 
for the case $\gamma(s)=\widehat{\zeta}(s)^2$, $\varepsilon=+1$ and $n=2$, 
which corresponds to the case of two-dimensional zeta integral over $\Q$, 
namely, the case of $h_{\cal E}(x)$, $Z_{\cal E}(x)$ and $\omega_{\cal E}(s)$. 
In particular, 
\begin{equation} \label{0301}
V_{\gamma}^{+1,2}(x,\nu)=V(x,\nu),
\end{equation}
where the right-hand side is defined in \eqref{007}. 
Recall that $Z_{\cal E}(x)$ involves the function  
\[
Z(x,\nu) = \Bigl(-x\frac{d}{dx} \Bigr)^4 V(x,\nu) 
= \Bigl(-x\frac{d}{dx} \Bigr)^4 V_{\gamma}^{+1,2}(x,\nu).  
\]

To state results simply, we use the non-holomorphic Eisenstein series. 
The (completed) non-holomorphic Eisenstein series $E^{\ast}(z, s)$ 
for the modular group ${\rm PSL}(2, \Bbb Z)$ 
is given for $z= x+iy$
with $y > 0$ and $\Re(s) > 1$ by 
\begin{equation}~\label{101}
E^{\ast}(z, s) = \frac{1}{2} \pi^{-s} \Gamma(s)  \!\!\!\!
\sum_{{(m,n) \in {\Bbb Z}^{2}}\atop{(m,n)\not=(0,0)}} \frac{y^s}{|mz+ n|^{2s}}. 
\end{equation}
It is well known that for fixed $z$, $E^{\ast}(z, s)$ is continued holomorphically  
to the whole $s$-plane except for two simple poles at $s=0$ and $s=1$
with residues $-1/2$ and $1/2$, respectively. 
As a function of $z$, $E^{\ast}(z, s)$ satisfies 
\begin{equation}~\label{102}
E^{\ast}\bigl(\, \frac{az+d}{cz+d}, s \bigr)= E^{\ast}(z, s) \quad 
\text{for all} \quad 
\gamma = 
\begin{pmatrix}
a & b \\ c & d
\end{pmatrix} 
\in {\rm PSL}(2,{\Bbb Z}).  
\end{equation}
We use notations 
\begin{equation}\label{104}
E(y)=E^\ast(iy,1/2) \quad \text{and} \quad  
Q=(4\pi)^{-1} e^{\gamma},
\end{equation} 
where $\gamma=0.57721+$ is the Euler's constant.
\begin{proposition}~\label{thm1}
Let $V(x,\nu)$ be the function defined by $\eqref{007}$. 
Then we have
\begin{align}\label{103}
V(x,\nu)
& = 
4 \sum_{N=1}^\infty \sigma_0(N) \bigl(  K_0(2 \pi N \nu x^{-2}) - x^2 K_0(2 \pi N \nu x^{2}) \bigr)
\\ \label{103_2} 
& = 
 x^2 \log x^2
+ x^2 \log Q \nu 
+ \log x^2
  - \log Q \nu  
+ \frac{x}{\sqrt{\nu}} 
\Bigl[
E( \nu x^{-2}) - 
E( \nu x^2)
\Bigr], 
\end{align}
where $K_0(t)$ is the $K$-Bessel function. 
\end{proposition}

As a corollary of Proposition \ref{thm1}, 
we obtain the asymptotic behavior of $Z(x, \nu)$ as $x \to 0^+$ and $+\infty$.

\begin{proposition}~\label{thm2}
We have
\begin{equation}\label{105}
\aligned
Z(x,\nu) & = 
16 x^2 \log x^2
+ 16 x^2 \log Qe^4 \nu
+ R_0(x,\nu), \\
R_0(x,\nu) & =
\frac{1}{2 \pi i} \int_{(c)} 
(2s)^4 \, \bigl[ (\nu^{-1} x^2)^s - \nu^{-1} (\nu x^2)^s \bigr] \, \widehat{\zeta}(s)^2 ds 
\quad (c>1)
\endaligned
\end{equation}
and 
\begin{equation}\label{106}
\aligned
Z(x,\nu)  & = 
\frac{16}{\nu} x^2 \log x^2 + \frac{16}{\nu}x^2 \log \frac{Q e^4}{\nu} 
+ R_\infty(x,\nu), \\
R_\infty(x,\nu) & =
\frac{1}{2 \pi i} \int_{(c)} 
(2s)^4 \, \bigl[ (\nu^{-1} x^2)^{s} - \nu^{-1} (\nu x^2)^{s} \bigr] \, \widehat{\zeta}(s)^2 ds 
\quad (c<0).
\endaligned
\end{equation}
Moving the path of integration to the line $\Re(s)=1+\delta$ or $\Re(s)=-\delta$ $(\delta>0)$, 
we have
\[
R_0(x,\nu) =  O\bigl( (\nu^{-1-\delta} + \nu^{\delta}) \, x^{2(1+\delta)} \bigr)
\]
for $0 < x \leq 1$, and 
\[
R_\infty(x,\nu) = O\bigl( (\nu^{-1-\delta} + \nu^{\delta}) \, {x^{-2\delta}} \bigr)
\]
for $x \geq 1$, respectively. 
Here the implied constants depend only on $\delta>0$. 
In particular, there exists $x_\nu>0$ and $x_\nu^\prime >0$ 
such that $Z(x,\nu)$ is negative for $0 < x <x_\nu$ and $Z(x,\nu)$ is positive for any $x > x_\nu^\prime$, 
respectively. 
\end{proposition}

An immediate consequence of Proposition \ref{thm2} is that for any $\nu^\prime \geq 1$
there exists $x_0=x_0(\nu^\prime) > 0$ such that 
\[
\sum_{\nu=1}^{\nu^\prime} c(\nu) Z(x,\nu) < 0 \quad x \in (0,x_0)
\]
if $c(\nu)>0$ for some $1 \leq \nu \leq \nu^\prime$ 
(recall that we assumed $c(\nu) \geq 0$). 
However, to decide the behavior of $Z_{\cal E}(x)$ near the zero, 
we need further considerations, 
since the infinite sum over the first term in the right-hand side of $\eqref{105}$ 
diverges for every $x>0$. 
Applying Proposition \ref{prop_004} to $Z_{\cal E}(x)$ we obtain the following. 

\begin{proposition} \label{prop_007}
Let ${\frak c}=\{c(\nu)\}$ be the nonnegative sequence defined by 
$\sum_{\nu=1}^\infty c(\nu) \nu^{-s} = c_{\cal E}^{1-2s} \zeta_{\cal E}(s)^2$, 
and $\gamma(s)=\widehat{\zeta}(s)^2$. In this case $\lambda_{\frak c}=\lambda_\gamma=2$. 
We denote $Z_{{\frak c},\gamma,i}^{2}(x)$ by $Z_{{\cal E},i}(x)$, 
where $Z_{{\frak c},\gamma,i}^n(x)~(i=0,1)$ are defined in \eqref{0208}. 
Then 
\begin{equation} \label{0309}
Z_{{\cal E},0}(x) = \frac{1}{2\pi i} \int_{(c)} 
\widehat{\zeta}(s)^2 \cdot c_{\cal E}^{1-2s} \zeta_{\cal E}(s)^2
\cdot (2-2s)^4 \cdot x^{2-2s} ds
\end{equation} 
and 
\begin{equation} \label{0310}
Z_{{\cal E},0}(x) = \frac{2}{\pi} 
\sum_{n=1}^\infty \frac{1}{n} 
\Bigl( \sum_{d|n} c(d)\sigma_{0}(n/d) \Bigr) 
{\cal K} (2 \pi n x^2)
\end{equation} 
where 
\begin{equation} \label{0311}
{\cal K} (x)=(16x^5+288x^3+16x)K_0(x)-(128x^4+64x^2)K_1(x). 
\end{equation}
\end{proposition}
\noindent
The series $Z_{{\cal E},0}(x)$ can be written as 
\[
Z_{{\cal E},0}(x)=\frac{2}{\pi}\sum_{n=1}^\infty {\cal K}(2 \pi n x^2)\frac{x^2(n+1) - x^2 n }{x^2n} \, a(n). 
\]
where $a(n)=\sum_{d|n} c(d)\sigma_{0}(n/d)$. 
This formula suggests that 
the limit $\displaystyle{\lim_{x \to 0^+} Z_{{\cal E},0}(x)}$ behaves like the Riemannian integral
\[
\lim_{x \to 0^+} \sum_{n=1}^\infty {\cal K}(2\pi n x^2) \frac{x^2(n+1) - x^2 n }{x^2n} 
= \int_{0}^{\infty} {\cal K}(2\pi x) \frac{dx}{x}. 
\]
We find that the integral on the right-hand side is zero, 
and $Z_{{\cal E},0}(x)$ tends to zero as $x \to 0^+$ 
if the truncated sum $\sum_{n \leq N} a(n)$ increases slowly. 
Hence the rough understanding of the limit behavior of $Z_{{\cal E},0}(x)$ may be considered as  
\[
\lim_{x \to 0^+} Z_{{\cal E},0}(x) = \int_{0}^{\infty} {\cal K}(2\pi x) \, d\mu_{\frak a} (x), 
\]
where the right-hand side is the ``integral'' of ${\cal K}(2\pi x)$ 
with respect to the ``measure'' $d\mu_{\frak a}$ attached to ${\frak a}=\{a(n)\}$. 
In fact, by using the partial summation, we have
\[
\sum_{n=1}^\infty \frac{a(n)}{n} {\cal K}(n x^2) 
 = -  \int_{0}^{\infty} 
\bigl( \sum_{n \leq v/x^2} \frac{a(n)}{n} \bigr)  {\cal K}^\prime(v)  dv 
 = \int_{0}^{\infty} 
 {\cal K}(2\pi u) d\mu_{{\frak a},x}(u),
\]
where 
$\mu_{{\frak a},x}(u) = \sum_{n \leq u/x^2} a(n)/n$.
\smallskip

The truncation version of \eqref{0310} can be used for a computational evidence of 
the single sign property (F-1). 

\begin{proposition} \label{prop008}
Let $T>0$ and $R>1$ be positive real numbers.  
Then, for any fixed $0<\varepsilon<1$ and $A>1$, we have 
\begin{equation}\label{501}
Z_{{\cal E},0}(x) = \frac{2}{\pi} 
\sum_{n \leq T} \frac{1}{n}
\Bigl( \sum_{d|n} c(d)\sigma_{0}(n/d) \Bigr) 
{\cal K} (2 \pi n x^2) +O ( T^{\varepsilon} (x^2 T)^{-A})
\end{equation}
for $x \geq \sqrt{R/T}$, where ${\cal K} (x)$ is in \eqref{0311} and the implied constant depends on $A$ and $\varepsilon$. 
In particular, for any $0<\alpha<1$ and $\beta>1$ such that $\alpha \beta > \varepsilon$,  we have
\begin{equation}\label{502}
Z_{{\cal E},0}(x) = \frac{2}{\pi} 
\sum_{n \leq R x^{-2-\alpha}} \frac{1}{n}
\Bigl( \sum_{d|n} c(d)\sigma_{0}(n/d) \Bigr) 
{\cal K}(2 \pi n x^2) + O ( x^{ \alpha \beta -\varepsilon } )
\end{equation}
for $0< x <1$, 
where the implied constant depends only on $\beta$ and $\varepsilon$. 
\end{proposition}
\begin{rem}
One suitable choice of $R$ for a numerical computation is $R=20$.
\end{rem}


\subsection{Proof of Proposition \ref{thm1}} 
At first we show \eqref{0301}. By definition \eqref{301_1}, 
\[
V_{\gamma}^{+1,2}(x) = \kappa_{\gamma}(\nu x^{-2}) - x^2 \kappa_{\gamma}(\nu x^2),
\]
where
\[
\kappa_\gamma(x) = \frac{1}{2\pi i} \int_{\Re(s)=c} \widehat{\zeta}(s)^2 x^{-s} ds \quad (c>1). 
\]
In general, if $F_i(s)=\int_{0}^{\infty}f_i(x)x^{s-1}dx~(i=1,2)$ for $\Re(s)>\delta$, then  
\[
\frac{1}{2\pi i} \int_{\Re(s)=c} F_1(s)F_2(s) x^{-s} ds 
= \int_{0}^{\infty} (f_1 \ast f_2)(x) \, x^{s}\frac{dx}{x} \quad (c>\delta)
\]
subject to appropriate conditions, where $(f_1 \ast f_2)(x)=\int_{0}^{\infty} f_1(a)f_2(xa^{-1})a^{-1}da$. 
Hence we have
\[
\kappa_{\gamma}(x) = 
\int_{0}^{\infty} 
\bigl( \theta(a^2) - 1 \bigr) \bigl( \theta(x^2a ^{-2} ) -1 \bigr) \frac{da}{a}, 
\]
since $\widehat{\zeta}(s)=\int_{0}^{\infty}(\theta(x^2)-1)x^{s-1}dx$ for $\Re(s)>1$. 
This shows
\[
\kappa_{\gamma}(\nu x^{-2}) - x^2 \kappa_{\gamma}(\nu x^2) 
= \int_{0}^{\infty} w_{a,\nu a^{-1}}(x) \frac{da}{a} = V(x,\nu)
\]
by definition \eqref{003} of $w_{a,\nu a^{-1}}(x)$, and we get \eqref{0301}. 
Using 
\[
\zeta(s)^2 = \sum_{N=1}^{\infty} \sigma_0(N) N^{-s} \quad (\Re(s)>1), 
\]
where $\sigma_\mu(n)=\sum_{d|n} d^\mu$, we have
\[
\kappa_\gamma(x) = \sum_{N=1}^{\infty} \sigma_0(N) 
\frac{1}{2 \pi i}\int_{(c)} \Gamma(s/2)^2  (\pi N x)^{-s} ds \quad (c>1).
\] 
The interchange of summation and integration is justified by Fubini's theorem. Using the formula  
$
\Gamma(s/2)^2 = 4 \int_{0}^{\infty} K_0(2x) x^{s-1} dx
$ (${\rm Re}(s)>0$)~\cite[p.\,14]{MR0350075} and the Mellin inversion formula, 
we have
\[
\frac{1}{2 \pi i} 
\int_{(c)} \Gamma(s/2)^2 x^{-s} ds 
= 4 K_{0}(2x)  \quad (c > 0).
\]
Hence we obtain
\begin{equation}\label{205}
\kappa_{\gamma}(x) = 4 \sum_{N=1}^{\infty} \sigma_0(N) K_0(2 \pi N x). 
\end{equation}
Combining this equality with $\eqref{0301}$, we obtain $\eqref{103}$. 
As for $\eqref{103_2}$, we refer the Fourier expansion of $E^\ast(z,s)$. 
For the particular case $z=iy$ and $s=1/2$, the Fourier expansion of $E^\ast(z,s)$ gives 
\begin{equation*}
E^\ast(iy,1/2) = \sqrt{y}\log y +(\gamma - \log 4 \pi)\sqrt{y} 
+ 4 \sqrt{y} \sum_{N=1}^\infty 
\sigma_{0}(N) K_{0}(2 \pi N y).
\end{equation*}
Combining this one with $\eqref{205}$, we obtain $\eqref{103_2}$. 
\hfill $\Box$


\subsection{Proof of Proposition \ref{thm2}} 
Recall equality $\eqref{103_2}$ in Proposition \ref{thm1} and notation $\eqref{104}$. 
We prove Proposition \ref{thm2} by using the modular relation $E(y)=E(y^{-1})$.  
Replacing $E(\nu x^2)$ by $E(\nu^{-1} x^{-2})$ in $\eqref{103_2}$, we have
\begin{equation}\label{213}
V(x,\nu) 
 = 
x^2 \log x^2
+ x^2 \log Q \nu
+ \log x^2
  - \log Q \nu  
+ \frac{x}{\sqrt{\nu}} 
\bigl[
E\bigl( \frac{\nu}{x^2} \bigr) - 
E\bigr( \frac{1}{\nu x^2} \bigr)
\bigr].
\end{equation}
Thus definition $\eqref{008}$ gives
\begin{equation}\label{214}
Z(x,\nu) 
 = 
16 x^2 \log x^2
+ 16 x^2 \log Qe^4 \nu
+ \bigl(x\frac{d}{dx}\bigr)^4  
\Bigl[ \frac{x}{\sqrt{n}} 
\bigl(
E\bigl( \frac{\nu}{x^2} \bigr) - 
E\bigr( \frac{1}{\nu x^2} \bigr)
\bigr) \Bigr].
\end{equation}
In the third term, we have
\begin{equation*}
\aligned
\frac{x}{\sqrt{\nu}} &
\bigl[  
E\bigl( \frac{\nu}{x^2} \bigr) - 
E\bigr( \frac{1}{\nu x^2} \bigr)
\bigr] \\
= & \, \bigl(1+\frac{1}{\nu}\bigr)\log \nu 
- \bigl(1-\frac{1}{\nu}\bigr)(\log x^2 - \log Q) 
 + 
\frac{1}{2 \pi i} \int_{(c)} 
\bigl[ (\nu^{-1} x^2)^s - \nu^{^-1} (\nu x^2)^s \bigr] \, \widehat{\zeta}(s)^2 ds
\endaligned
\end{equation*}
by using the Fourier expansion of the Eisenstein series and \eqref{205}. 
Therefore
\begin{equation}\label{215}
\aligned
\bigl(x\frac{d}{dx}\bigr)^4  
\Bigl( \frac{x}{\sqrt{\nu}} 
\bigl[ 
E\bigl( & \frac{\nu}{x^2} \bigr) - 
E\bigr( \frac{1}{\nu x^2} \bigr)
\bigr] \Bigr) \\
& = \frac{1}{2 \pi i} \int_{(c)} 
\bigl(x\frac{d}{dx}\bigr)^4\bigl[ (\nu^{-1} x^2)^s - \nu^{-1} (\nu x^2)^s \bigr] \, \widehat{\zeta}(s)^2 ds \\
& = \frac{1}{2 \pi i} \int_{(c)} 
(2s)^4 \, \bigl[ (\nu^{-1} x^2)^s - \nu^{-1} (\nu x^2)^s \bigr] \, \widehat{\zeta}(s)^2 ds
\endaligned
\end{equation}
for $c>1$. Hence we obtain \eqref{105}. 
Formula \eqref{106} is obtained by a way similar to the above arguments. 
In that case, we replace $E(\nu x^{-2})$ by $E(\nu^{-1} x^2)$ in $\eqref{103_2}$.  
\hfill $\Box$ 

\subsection{Proof of Proposition \ref{prop_007}} 
At first we note that $\lambda_{\frak c}=\lambda_\gamma=2$ in this case. 
The first equality is a direct consequence of definition \eqref{0207} and \eqref{0208}. 
By \eqref{0207} and \eqref{205} we have
\[
\aligned
Z_{{\cal E},0}(x) 
& = 
\Bigl( -x\frac{d}{dx} \Bigr)^4 
\left( 4 x^2 \sum_{n=1}^{\infty} \bigl( \sum_{d|n} c(d) \, \sigma_{0}(n/d) \bigr)  K_{0}(2\pi n x^2) \right) \\
& = \frac{2}{\pi}
\Bigl( -x\frac{d}{dx} \Bigr)^4 
\left(\sum_{n=1}^{\infty} \frac{1}{n} \bigl( \sum_{d|n} c(d) \, \sigma_{0}(n/d) \bigr) (2\pi n x^2 K_{0}(2\pi n x^2)) \right).
\endaligned
\]
Hence the proof of the second equality will complete by the following lemma. 
\begin{lemma} Let $A$ be a positive real number. 
Then we have
\begin{equation}\label{211}
\bigl(x\frac{d}{dx}\bigr)^4 \bigl[ K_0(Ax^{-2}) \bigr]
=   \bigl(\frac{64A^2}{x^4} + \frac{16 A^4}{x^8}\bigr) K_0(Ax^{-2})
  - \frac{64 A^3}{x^6}K_1(Ax^{-2})
\end{equation}
and 
\begin{equation}\label{212}
\aligned
\bigl(x\frac{d}{dx}\bigr)^4  \bigl[ A x^2 K_0(Ax^2) \bigr]
= & \, (16 A^5 x^{10} + 288 A^3 x^6 + 16 A x^2) K_0(Ax^{2})  \\
  & - (128 A^4 x^8 + 64 A ^2x^4) K_1(Ax^{2}) .
\endaligned
\end{equation}
\end{lemma}
\begin{pf}
Using the formula 
\[
\frac{d}{dz} [z^\nu K_\nu(z)] = -z^\nu K_{\nu-1}(z), \quad
\frac{d}{dz} [z^{-\nu} K_\nu(z)] = -z^{-\nu} K_{\nu+1}(z)
\]
and $K_0^\prime(z)=-K_1(z)$, we have
\[
\bigl(x\frac{d}{dx}\bigr)^4 \bigl[ K_0(x^{-2}) \bigr]
=   \bigl(\frac{64}{x^4} + \frac{16}{x^8}\bigr) K_0(x^{-2})
  - \frac{64}{x^6}K_{1}(x^{-2}).
\]
This implies $\eqref{211}$, 
since the multiplication $f(x) \mapsto f(A x)$ 
and the differential $x\frac{d}{dx}$ are commutative each other.  
By a way similar to this, we obtain $\eqref{212}$.
\end{pf}

\subsection{Proof of Proposition \ref{prop008}}
We denote $\sum_{d|n} c(d)\sigma_{0}(n/d)$ by $a(n)$. 
Then, for any $\varepsilon>0$, there exists $M_{\varepsilon}>0$ 
such that $|a(n)| \leq M_{\varepsilon} n^{\varepsilon}$
The asymptotic formula of $K_{\nu}(t)$ for $t>0$ and fixed $\nu \in \C$ is given by
\begin{equation}\label{210}
K_\nu(t) = \Bigl( \frac{\pi}{2t} \Bigr)^{1/2} e^{-t} \bigl( 1+ \frac{\theta}{2t} \bigr),
\end{equation}
where $|\theta| \leq |\nu^2 - (1/4)|$ (see (23.451.6) of~\cite{MR2360010}). 
Using $\eqref{210}$ we have
\[
\aligned
\sum_{n > T} \frac{a(n)}{n} {\cal K}(x^2 n) 
& \leq  \sqrt{\frac{\pi}{2}} \sum_{n > T} \frac{a(n)}{n} 
 (16 x^{10} n^5 + 288 x^6 n^3 + 16 x^2 n) \frac{e^{-x^2 n}}{\sqrt{x^2 n}} \bigl(1+\frac{|\theta_0|}{2 x^2 n}\bigr) \\
& + \sqrt{\frac{\pi}{2}} \sum_{n > T} \frac{a(n)}{n}
 (128 x^8 n^4 + 64 x^4 n^2) \frac{e^{- x^2 n}}{\sqrt{x^2 n}} \bigl(1+\frac{|\theta_1|}{2 x^2 n}\bigr). 
\endaligned
\]
Also
\[
\sum_{n > T} \frac{a(n)}{n} {\cal K}(x^2 n) 
 \leq 312 \sqrt{2 \pi} \, x ^9 \sum_{n \geq T} a(n) n ^{7/2}  e^{-x^2 n},
\]
since $|\theta_0| \leq 1/4$, $|\theta_1| \leq 3/4$ and $ x^2 n \geq x^2 T \geq R > 1$.
Using the estimate $  t^{-k-(9/2)} \geq t^{-k-5} \geq e^{-t}/(k+5)!  $ for $t > 0$, 
we have 
\[
x ^9 \sum_{n > T} a(n) n ^{7/2}  e^{-x^2 n} 
\leq (k+5)! \, x^{-2 k} \sum_{n \geq T} a(n) n^{-1-k}. 
\]
Therefore
\[
\sum_{n > T} \frac{a(n)}{n} {\cal K}(x^2 n)  \leq 
312 \sqrt{2 \pi} \, (k+5)! \, x^{-2 k} \sum_{n \geq T} a(n) n^{-1-k}.
\]
Because of $|a(n)| \leq M_\varepsilon n^\varepsilon$, we obtain
\[
\sum_{n > T} a(n) n^{-1-k} 
\leq M_\varepsilon  \sum_{n >T} n^{-1- k + \varepsilon} 
\leq \frac{M_\varepsilon}{k-\varepsilon} T^{\varepsilon -k}.
\]
Together with the above, we have
\[
x^{-2 k} \sum_{n \geq T} a(n) n^{-1-k} \leq 
312 \sqrt{2 \pi} \, (k+5)! \, \frac{M_\varepsilon}{k-\varepsilon} T^{\varepsilon} (x^2 T)^{-k}. 
\]
This inequality gives $\eqref{501}$. 
\hfill $\Box$ 


\section{Proof of Theorem \ref{thm_001}} 
In this part, we prove the following Thorem \ref{prop4} which include Theorem \ref{thm_001}. 

\begin{theorem}\label{prop4}
Let $E$ be an elliptic curve of conductor $q_E$ over ${\Q}$, 
and let $\{c(\nu)\}$ be the nonnegative sequence defined by 
\[
\sum_{\nu=1}^\infty c(\nu) \nu^{-s} 
= c_{\cal E}^{1-2s} \zeta_{\cal E}(2s)^2
= c_{\cal E}^{1-2s} n_{\cal E}(2s)^2 \zeta_{E}(2s)^2, 
\]
where $\zeta_E(s)$ and $n_{\cal E}(s)$ are given by \eqref{001} and \eqref{002_1}, respectively.  
Let 
\[
Z_{\cal E}(x) = \sum_{\nu=1}^{\infty} c(\nu) Z(x,\nu).
\]
Suppose that the Riemann hypothesis for $L(E,s)$, 
namely, all nontrivial zeros of $L(E,s)$ lie on the line $\Re(s)=1$. 
Then we have 
\[
Z_{\cal E}(x) = O(x^{1-\varepsilon})
\]
for any $\varepsilon>0$, where implied constant depends on $E$ and $\varepsilon$.  

Further, suppose that all zeros of $L({\cal E},s):=n_{\cal E}(s)^{-1}L(E,s)$ 
are simple except for the possible zero at $s=1$ and the estimate 
\begin{equation}\label{606}
\sum_{0< \Im(\rho) \leq T}|L^\prime({\cal E},\rho)|^{-2} = O(T)
\end{equation}
holds, where $\rho$ runs all zeros of $L({\cal E},s)$ on the line $\Re(s)=1$. 
Then we have
\begin{equation}\label{606_2}
Z_{\cal E}(x) = 
-C \, x ( \log (1/x))^{2r+2J+1}(1 + O((\log (1/x))^{-1}), 
\end{equation}
where $r$ is the order of $L(E,s)$ at $s=1$, $J (\geq 1)$ is the number of Euler factors in $n_{\cal E}(s)$ 
and $C$ is a positive constant. 
In particular, there exists $x_{\cal E} >0$ 
such that $Z_{\cal E}(x)$ is negative on $(0,x_{\cal E})$. 
\end{theorem}

\begin{theorem}\label{prop04}
Let $E$ be an elliptic curve of conductor $q_E$ over ${\Q}$, 
and let $\{c(\nu)\}$ be the nonnegative sequence defined by 
$\sum_{\nu=1}^\infty c(\nu) \nu^{-s} = q_E^{-2s}\zeta_{E}(2s)^2$, 
where $\zeta_E(s)$ is given by \eqref{001}. 
Denote by $r$ the order of $L(E,s)$ at $s=1$.  
Let $Z_{E}(x) = \sum_{\nu=1}^{\infty} c(\nu) Z(x,\nu)$.
Suppose that the Riemann hypothesis for $L(E,s)$, 
and all zeros of $L(E,s)$ are simple 
except for the possible zero at $s=1$ and the estimate \eqref{606} holds for $L^\prime(E,s)$. 
Then 
\begin{equation}
Z_E(x) = 
\begin{cases}
-C \, x ( \log (1/x))^{2r+1}(1 + O((\log (1/x))^{-1}) & \text{if $r \geq 1$}, \\
-(C+v(x)) \, x \log (1/x) \, (1 + O((\log (1/x))^{-1}) & \text{if $r=0$}, 
\end{cases}
\end{equation}
where $C$ is a positive constant and $v(x)$ is a bounded function. 
In particular, if $L(E,1)=0$, there exists $x_E >0$ 
such that $Z_{E}(x)$ is negative on $(0,x_{E})$. 
\end{theorem}

\begin{rem}
Computational table of $Z_E(x)$ and its {\rm PARI/GP} program can be available from 
Fesenko's homepage~\cite{Fe4} 
(for convenience, use it together with the table of elliptic curves in Appendix B.5 of Cohen~\cite{MR1228206}). 
These numerical table suggests that $Z_E(x)$ is negative for sufficiently small $x>0$ 
even in the case that the $\Z$-rank of the Mordell-Weil group $E({\Q})$ is zero. 
\end{rem} 

\begin{rem}
Consequeces 
similar to Theorem \ref{prop4} and Theorem \ref{prop04} 
are obtained for the case of the general number field 
under conditions analogous to Theorem \ref{prop4} and Theorem \ref{prop04} 
if we add the assumption that 
$L(E,s)$ is extended to an entire function 
and satisfies the conjectural functional equation. 
The proof goes very similarly 
using standard analytic properties of $L$-functions (cf. \cite[section 5]{IK}). 
\end{rem}
Estimate \eqref{606} is the analogue of the conjectural estimate
\begin{equation}\label{607}
\sum_{0 < \gamma \leq T} |\zeta^\prime(\rho)|^{-2k} = O(T(\log T)^{(k-1)^2}) 
\quad (k \in {\Bbb R}),
\end{equation}
where we assume that the Riemann hypothesis and all zeros of $\zeta(s)$ are simple. 
Estimate \eqref{607} was independently conjectured by Gonek~\cite{MR1014202} and Hejhal~\cite{MR993326} 
from different points of view. For $k=1$ Gonek conjectuted the asymptotic formula
$\sum_{0 < \gamma \leq T}|\zeta^\prime(\rho)|^{-2} \sim (3/\pi^3) \, T$.  
Let $f$ be a normalized Hecke eigenform of weight $k>1$ and level $q$ with trivial nebentypus. 
Murty and Perelli~\cite{MR1692847} had shown that almost all zeros of $L(f,s)$ are simple 
if we assume the Riemann hypothesis for $L(f,s)$ and the pair correlation conjecture for it. 
According to the Shimura-Taniyama-Weil conjecture proved by Wiles et al., 
almost all zeros of $L(E,s)$ are simple 
if we assume the Riemann hypothesis of $L(E,s)$ and its pair correlation conjecture.  
\smallskip

We prove Theorem \ref{prop4} only, 
since Theorem \ref{prop04} is proved by a similar way. 
To prove Theorem \ref{prop4}, we need the following lemma. 

\begin{lemma}\label{lem601}
Let $E$ be an elliptic curve over ${\Q}$. Let $H \geq 1$ be a real number. 
Then there exists a real number $A$ and a subset ${\frak S}_T$ of $[\,T,T+1)$ such that 
\[
|L(E,\sigma \pm it)|^{-1} =O(t^{A}) \quad (-1/2 \leq \sigma \leq 5/2, ~ t \in {\frak S}_T)
\]  
and the Lebesgue measure of $[\,T,T+1) \setminus {\frak S}_T$ is less than or equal to $1/H$.  
\end{lemma}
\begin{pf}
This is shown by a way similar to the proof of Theorem 9.7 in Titchmarsh~\cite{MR882550},  
because of the modularity of $E/{\Q}$ (see also \cite[section 5.1]{IK}).    
\end{pf}
Under the Riemann hypothesis for $L(E,s)$ 
we can obtain more sharp estimate for $L(E,s)^{-1}$ 
in a vertical strip, but we do not need such sharp estimate for Theorem \ref{prop4}. 
\smallskip

\noindent 
{\bf Proof of Theorem \ref{prop4}.} 
It suffices to deal with $Z_{{\cal E},0}(x)$ in Proposition \ref{prop_007}, 
since 
\begin{equation}\label{609}
Z_{\cal E}(x) = - Z_{{\cal E},0}(x) + O(x^A) \quad (x \to 0^+)
\end{equation} 
by Proposition \ref{prop_004}. 
We denote by $\rho=1+i\gamma$ the zeros of $L({\cal E},s)=n_{\cal E}(s)^{-1}L(E,s)$ 
in $0 < \Re(s) < 2$ (on $\Re(s)=1$). 
At first we prove that each interval $[n,n+1)$ 
contains a value $T$ for which 
\begin{equation}\label{610}
\aligned
Z_{{\cal E},0}(x) 
= & \,\,
x \, P_1(\log (1/x))  + x^2 ( C_2 \log (1/x) + C_2^\prime) \\
& + \sum_{0<|\gamma| \leq T} x^{1-i\gamma} ( C_\gamma \log (1/x) + C_\gamma^\prime) \\
& + O(x^{-\delta}e^{-(\pi/4-\delta)T}) + O(x^{2+\delta}),   
\endaligned
\end{equation}
for any fixed $0<\delta<\pi/4$, where $P_1$ is a polynomial of degree $2r+2J+1$, 
$C_\gamma$, $C_\gamma^\prime$, $C_2$ and $C_2^\prime$ are constants. 
Taking $T=1/x$ in \eqref{610} and then tending $x$ to zero, we obtain
\begin{equation}\label{610_2}
\aligned
Z_{{\cal E},0}(x) 
= & \,\,
x \, P_1(\log (1/x))  
+ O(x(\log (1/x))) \quad (x \to 0^+).
\endaligned
\end{equation}
This asymptotic formula shows \eqref{606_2} via \eqref{609} except for the sign of $C$ in \eqref{606_2}. 
The positivity of $C$ follows from \eqref{leading} in below.  

Now we prove \eqref{610}. 
We divide the ingetral $\int_{(c)}$ in \eqref{0309} into three parts 
$\int_{c-iT}^{c+iT}$, $\int_{c+iT}^{c+i\infty}$ and $\int_{c-i\infty}^{c-iT}$. 
We consider the positively oriented rectangle 
with vertices at $c+iT$, $c^\prime +iT$, $c^\prime -iT$ and $c-iT$ with $-1<c^\prime<0$. 
In this rectangle the integrand has poles 
$s=1$ of order $2r+2J+2$, 
$s=0$ of order $2$, 
and $s=1+i\gamma$ of order $2$, 
and has no other poles, 
since we assumed that all zeros of $L({\cal E},s)$ are simple except for $s=1$. 
Thus the residue theorem gives
\begin{equation}\label{611}
\aligned
\frac{1}{2\pi i}\int_{c-iT}^{c+iT} ds
= & \,\, x \, P_1(\log (1/x))  
+  x^2 ( C_2 \log (1/x) + C_2^\prime)\\
& +  \sum_{0<|\gamma| \leq T} x^{1-i\gamma} ( C_\gamma \log (1/x) + C_\gamma^\prime)\\ 
& + \frac{1}{2\pi i} \int_{c^\prime-iT}^{c^\prime+iT} ds 
+ \frac{1}{2\pi i} \int_{c^\prime +iT}^{c+iT} ds 
- \frac{1}{2\pi i} \int_{c^\prime -iT}^{c-iT} ds,
\endaligned
\end{equation}
where $P_1$ is the polynomial of degree $2r+2J+1$ given by
\[
P_1(\log (1/x)) 
= x^{-1} \underset{s=1}{\rm Res} 
\Bigl(\, 
\widehat{\zeta}(s/2)^2 c_{\cal E}^{1-s} \zeta_{\cal E}(s)^2 (s-2)^4 x^{2-s} \,\Bigr) \in {\R}[\log(1/x)]
\]
and 
\[
C_\gamma \log (1/x) + C_\gamma^\prime 
= x^{-1+i\gamma}\underset{s=1+i \gamma}{\rm Res} 
\Bigl(\, 
\widehat{\zeta}(s/2)^2 c_{\cal E}^{1-s} \zeta_{\cal E}(s)^2 (s-2)^4 x^{2-s} \,\Bigr).
\]
In particular the leading term of $P_1$ is given by
\begin{equation}\label{leading}
\widehat{\zeta}(1/2)^2 \zeta(0)^2 \lim_{s \to 1}\Bigl[ \frac{L({\cal E},s)}{(s-1)^{r+J}} \Bigr]^{-2} 
(\log (1/x))^{2r+2J+1}.
\end{equation}
To calculate the integrals in the right-hand side of \eqref{611}, 
we use Lemma \ref{lem601}, the well-known (unconditional) estimate
\[
\zeta(\sigma+it) \ll 
\begin{cases}
\,1 & \sigma>1, \\ 
|t|^{(1-\sigma)/2+\varepsilon} & 0 \leq \sigma \leq 1, ~|t| \geq 2, \\
|t|^{1/2-\sigma} & \sigma<0, ~|t| \geq 2, \\
\end{cases}
\]
and Stirling's formula
\[
\Gamma(s) = \sqrt{2\pi} \, |t|^{\sigma-1/2}e^{-(\pi/2)|t|}(1+O(|t|^{-1}))
\quad (\sigma_1 \leq \sigma \leq \sigma_2,~|t| \geq 1).
\] 
Using these estimates the second and the third integrals in the right-hand side of \eqref{611} are 
\[
\frac{1}{2\pi i} \int_{c^\prime +iT}^{c+iT} ds 
- \frac{1}{2\pi i} \int_{c^\prime -iT}^{c-iT} ds
\ll x^{2-\sigma} |t|^{2(6+A)} e^{-(\pi/4)|t|}(1+O(|t|^{-1})) 
\ll x^{2-c} e^{-(\pi/4-\delta)T}
\]
for $c^\prime \leq \sigma \leq c$, $|t| \geq t_\delta >1$. 
Using the functional equation of $\zeta(s)$ and $\zeta_{\cal E}(s)$ 
(\cite[47 in \S4.2]{Fe3}), 
the first integral in the right-hand side of \eqref{611} is calculated as
\[
\aligned
\frac{1}{2\pi i} \int_{c^\prime-iT}^{c^\prime+iT} ds 
& = \frac{1}{2\pi} \int_{2-c^\prime-iT}^{2-c^\prime+iT}
\widehat{\zeta}(s/2)^2 c_{\cal E}^{1-s}\zeta_{\cal E}(s)^2 s^4
x^{s}ds \\
\ll & \,\, x^{2-c^\prime} \int_{2-c^\prime-i\infty}^{2-c^\prime+i\infty} 
 |\Gamma(s/4)|^2 |\zeta(s/2)\zeta_{\cal E}(s)|^2 |s|^4|ds| \ll x^{2-c^\prime}. 
\endaligned
\]
The last inequality follows from the fact that 
the Dirichlet series of $\zeta(s/2)\zeta_{\cal E}(s)$ converges absolutely for $\Re(s)>2$, 
Lemma \ref{lem601} and the above estimates for $\zeta(s)$ and $\Gamma(s)$. 
Finally we have 
\[
\frac{1}{2\pi i} \left( \int_{c+iT}^{c+i\infty} + \int_{c-i\infty}^{c-iT} \right)
\widehat{\zeta}(s/2)^2 c_{\cal E}^{1-s}\zeta_{\cal E}(s)^2 (s-2)^4 x^{2-s} ds 
\ll x^{2-c} e^{-(\pi/4 -\delta)T}.
\]
Combining the above three estimates and taking $c=2+\delta$ and $c^\prime=-\delta$, 
we obtain formula \eqref{610}. 

To justify \eqref{610_2}, we estimate the sum in the right-hand side of \eqref{610}. 
By an elementary calculation we obtain
\[
C_\gamma 
= \frac{1}{L^\prime({\cal E},\rho)^2} f(\rho)
\quad \text{and} \quad 
C_\gamma^\prime 
= \frac{L^{\prime\prime}({\cal E},\rho)}{L^\prime({\cal E},\rho)^3}f(\rho)
+ \frac{1}{L^\prime({\cal E},\rho)^2}f^\prime(\rho),
\]
where $f(s)= c_{\cal E}^{1-s} \widehat{\zeta}( s/2 )^2 \zeta( s )^2 \zeta( s-1 )^2 (s -2)^4$. 
Functions $f(s)$ and $f^\prime(s)$ are bounded by $e^{-A|t|}$ for some $A>0$ 
on the vertical line $s=1+it$. Therefore
\[
\aligned
\sum_{0<|\gamma| \leq T} |C_\gamma| 
& \ll \sum_{0<|\gamma| \leq T} |L^\prime({\cal E},\rho)|^{-2} e^{-A_1|\gamma|}
\ll \int_{1}^{T} \Bigl( \sum_{0<|\gamma| \leq t} |L^\prime({\cal E},\rho)|^{-2} \Bigr) e^{-A_1t} dt.
\endaligned
\]
Using assumption \eqref{606} in the right-hand side, we have
\begin{equation}\label{613}
\sum_{0<|\gamma|\leq T} |C_\gamma| \ll 1 + e^{-A_2 T}.
\end{equation}
On the other hand we have $|L^{\prime}({\cal E},\rho)|^{-1} = O(T^{1/2})$ for $0<|\gamma| \leq T$, 
because 
$|L^{\prime}({\cal E},\rho)|^{-2} 
\leq \sum_{0< |\gamma| \leq T} |L^{\prime}({\cal E},\rho)|^{-2} =O(T)$ 
by \eqref{606}. 
Therefore
\[
\sum_{0< |\gamma| \leq T} |L^{\prime}({\cal E},\rho)|^{-3} 
\ll T^{1/2} \sum_{0< |\gamma| \leq T} |L^{\prime}({\cal E},\rho)|^{-2} 
\ll T^{3/2}. 
\]
While we have the rough estimate $L^{\prime\prime}({\cal E},\rho)=O(|\gamma|^{3/2})$ 
by using the Cauchy estimate for derivatives and the Phragm\'en-Lindel\"of principle. 
Hence we have
\[
\sum_{0< |\gamma| \leq T} \Bigl| \frac{L^{\prime\prime}({\cal E},\rho)}{L^\prime({\cal E},\rho)^3} \Bigr| =O(T^3).
\]
Using this estimate we obtain
\begin{equation}\label{614}
\sum_{0<|\gamma|\leq T} |C_\gamma^\prime| \ll 1 + e^{-A_3T}
\end{equation}
by a way similar to the proof of \eqref{613}. 
By \eqref{613} and \eqref{614} we obtain
\[
\sum_{0<|\gamma| \leq T} x^{1-i\gamma} ( C_\gamma \log x + C_\gamma^\prime)
\ll x \log x \, (1 + e^{-A_4 T}).
\]
This estimate justify the process leading the asymptotic formula \eqref{610_2} from \eqref{610}. 

The first assertion in the proposition is obtained from \eqref{0309} 
by moving the path of integration to the vertical line $\Re(s)=1+\varepsilon$. 
It is justified by Lemma \ref{lem601} and the above estimates for $\zeta(s)$ and $\Gamma(s)$, 
and the resulting integral is estimated as $O(x^{1-\varepsilon})$ by the same tools. 
\hfill $\Box$

\section{On single sign property of $Z(x)$} 

Finally, we remark on the single sign property of $Z_{{\frak c},\gamma}^{+1,2}(x)$ with $\gamma(s)=\widehat{\zeta}(s)^2$ 
from a viewpoint of Euler products. In this section, 
we denote $Z_{{\frak c},\gamma}^{+1,2}(x)$ simply as $Z_{\frak c}(x)$ or $Z(x)$ 
if there is no confusion. 

\subsection{Euler product of degree 2}
Let $D=\{ s \in {\C} ~|~ \Re(s) > 1/2\}$, 
and denote by $H(D)$ the space of holomorphic functions on $D$ 
equipped with the topology of uniform convergence on compacta. 
Let $S$ be a finite set of primes. 
Let $\gamma=\{s\in {\C}~|~|s|=1\}$, 
and let 
\[
\Omega_S = \prod_{p \not\in S} \gamma_p,
\]
where $\gamma_p=\gamma$ for all primes $p \not \in S$. 
The infinite dimensional torus $\Omega_S$ is a compact topological Abelian group. 
Denote by $m_H$ the probability Haar measure on $(\Omega_S, {\cal B}(\Omega_S))$, 
where ${\cal B}(\Omega_S)$ is the class of Borel sets of the space $\Omega_S$. 
Thus we obtain probability space $(\Omega_S, {\cal B}(\Omega_S),m_H)$. 
Let $\omega(p)$ be the projection of $\omega \in \Omega_S$ to the coordinate space $\gamma_p$. 
For $\omega \in \Omega_S$, 
we define $L_S(s,\omega)$ by the Euler product 
\begin{equation}\label{904}
L_S(s,\omega) 
= \prod_{p \not\in S} \left( 1 - (\omega(p)+\overline{\omega(p)})p^{-s} + p^{-2s} \right)^{-1}. 
\end{equation}
The right-hand side converges absolutely for $\Re(s)>1$, 
since $|\omega(p)| = 1$. 
Thus $L_S(s,\omega)$ is a holomorphic function on $\Re(s)>1$ and has no zero in there. 
Moreover, for almost all $\omega \in \Omega_S$ with respect to $m_H$, the sum
\begin{equation}\label{905}
\sum_p \omega(p)p^{-s}
\end{equation}
converges for $\Re(s)>1/2$ (cf.~\cite[\S 5.5.1]{MR1376140}). 
In particular, for almost all $\omega \in \Omega_S$, 
\[
\log L_S(s,\omega) = - \sum_{p} \log\left( 1 - (\omega(p)+\overline{\omega(p)})p^{-s} + p^{-2s} \right) \in H(D)
\]
and $L_S(s,\omega) \in H(D)$ without zeros. 

We consider ${\frak c}_\omega=\{c_\omega(\nu)\}$ defined by 
\[
D_\omega(s)=\sum_{\nu=1}^{\infty} c_{\omega}(\nu) \, \nu^{-s} = \frac{\zeta(2s)^2\zeta(2s-1)^2}{L(2s-1/2,\omega)^2}. 
\]
The sequence ${\frak c}_\omega$ is nonnegative for all $\omega  \in \Omega_S$, 
since 
\[
\frac{1}{(1-p^{-s})(1-p^{1-s})} = \sum_{n=0}^{\infty} \left(\, \sum_{k=0}^{n} p^k \,\right)p^{-ns} 
\]
for $p \in S$, 
\[
\frac{1-(\omega(p)+\overline{\omega(p)})p^{1/2}p^{-s} + p^{1-2s}}{(1-p^{-s})(1-p^{1-s})}
= 1 + \left(\, p+1- 2 \sqrt{p}~ \Re (\omega(p)) \right) \sum_{n=1}^{\infty} \left(\, \sum_{k=0}^{n-1} p^k \right) p^{-ns}
\]
for $p \not\in S$, and $|\Re(\omega(p))| \leq 1$. 

As mentioned above, $L(s,\omega)$ is continued holomorphically 
to $\Re(s)>1/2$ without zeros for almost all $\omega \in \Omega_S$, 
because of the convergence of the series $\eqref{905}$. 
Hence,  we have
\begin{equation}\label{906}
Z_{\frak c} (x) = O(x^{1-\varepsilon}) \quad \text{and} \quad 
Z_{\frak c} (x) = \Omega(x) \quad \text{as} \quad x \to +0
\end{equation}
for almost all $\omega \in \Omega_S$. 
We denote by $A$ the set of all $\omega \in \Omega_S$ such that 
the series $\eqref{905}$ converges for $\Re(s)>1/2$. 
Also, we denote by $A^\prime$ the set of all $\omega \in \Omega_S$ such that 
the series $L_S(s,\omega)$ is continued holomolphically to $\Re(s)>1/2$ without zeros. 
Clearly, $A^\prime \supset A$, so $m_H(\Omega_S \setminus A^\prime) = m_H(\Omega_S \setminus A) = 0$. 

On the other hand, it is known that 
the vertical line $\Re(s)=1/2$ is a natural boundary 
of meromorphic continuation of $L_S(s,\omega)$ for almost all $\omega \in \Omega_S$. 
Therefore, for almost all $\omega \in \Omega_S$, 
the corresponding $Z_{\frak c}(x)$ may have oscillation near zero, 
even if the estimate $\eqref{906}$ holds. 
We denote by $B$ the set of all $\omega \in \Omega_S$ such that 
the line $\Re(s)=1/2$ is a natural boundary 
of meromorphic continuation of $L_S(s,\omega)$. 

Unfortunately, the relation among the sets $A$, $A^\prime$ and $B$ is not clear in general. 
However, the above considerations suggest that 
the single sign property of $Z_{{\frak c}_{\omega}}(x)$ is a special phenomenon and has a proper reason. 
Of course there is a possibility that $Z_{{\frak c}_{\omega}}(x)$ has the single sign property accidentally. 
To obtain a plausible reason for the single sign property, 
we should deal with some regular class of coefficients. 

\subsection{A conjecture related to the Selberg class}

The Selberg class $\cal S$ is a general class of Dirichlet series 
satisfying five axiom~\cite{MR1220477,MR2175035}. 
Roughly speaking, the Selberg class is a class of Dirichlet series 
having an Euler product, analytic continuation and functional equation of certain type. 
All known examples of functions in the Selberg class 
are automorphic (or at least conjecturally automorphic) $L$-functions. 
It is expected that all $L$-functions in the Selberg class 
satisfy an analogue of the Riemann hypothesis. 

For the Dirichlet series $D(s)$ equipped with the Euler product
\begin{equation}\label{907}
D(s) = \sum_{\nu=1}^{\infty} \frac{c(\nu)}{\nu^s} 
= \prod_{p} \left( 1 + \sum_{k=1}^\infty \frac{c(p^k)}{p^{ks}} \right),
\end{equation}
we take
\begin{equation}\label{908}
L(s)= \frac{\zeta(s+1/2) \zeta(s-1/2)}{ D(s/2+1/4)^{1/2}}. 
\end{equation}
Using the Euler product of $D(s)$ and formula $(1+X)^{-1/2} = 1 -\frac{1}{2} X + \frac{3}{4}X^2 - \cdots$, 
we obtain the Dirchlet series expansion of $L(s)$ and its Euler product
\begin{equation}\label{909}
L(s) = \sum_{n=1}^{\infty} \frac{\omega(n)}{n^s} 
= \prod_{p} \left( 1 + \sum_{k=1}^\infty \frac{\omega(p^k)}{p^{ks}} \right).
\end{equation}
Conversely, if we take
\begin{equation}\label{910}
D(s) = \left( \frac{\zeta(2s) \zeta(2s-1)}{ L(2s-1/2)^2} \right)^2,
\end{equation}
for the Dirichlet series $L(s)$ having an Euler product, 
then $D(s)$ is also a Dirichlet series equipped with the Euler product. 

Suppose that $D(s)$ and $L(s)$ are related as $\eqref{908}$ and $\eqref{910}$, 
and $D(s)$ or $L(s)$ has an Euler product.  
If $D(s)$ or $L(s)$ is continued meromorphically to a region containing $\Re(s) \geq 1/2$, 
then properties    
\begin{enumerate}
\item[(D1)] $D(s)$ has double pole at $s=1$,
\item[(D2)] $D(s)$ has the pole of order $2m+2$ at $s=1/2$,
\item[(D3)] $D(s)$ has the pole of order $2n$ at $s=\sigma \in (1/2,1)$,
\end{enumerate}
and 
\begin{enumerate}
\item[(L1)] $L(s)$ is regular at $s=1$,
\item[(L2)] $L(s)$ has the zero of order $m$ at $s=1/2$,
\item[(L3)] $L(s)$ has the zero of order $n$ at $s=\sigma \in (1/2,1)$,
\end{enumerate}
are equivalent to each other. 
If $Z(x)$ attached to $D(s)$ has a single sign for sufficiently small $x>0$, 
then $D(s/2)$ has the real pole at the abscissa $\sigma_c<2$ of convergence of $\eqref{305}$, 
and has no pole in the strip $\sigma_c<\Re(s)<2$. 
Under the Riemann hypothesis for $\zeta(s)$, 
this implies that $L(s)$ has no zero in the strip $\sigma_c/2<\Re(s)<1$. 
Conversely, if $L(s)$ has zero at $1/2 \leq \sigma_0<1$ and has no zero in the strip $\sigma_0<\Re(s)<1$, 
then $D(s/2)$ has pole at $1 \leq 2\sigma_0<2$ and has no pole in the strip $2\sigma_0<\Re(s)<2$. 
This suggests the single sign property of $Z(x)$ attached to $D(s)$. 
Through such relations, we interpret the expectation that all $L(s) \in {\cal S}$ satisfy the Riemann hypothesis 
in the language of $D(s)$ via $\eqref{908}$ and $\eqref{910}$.  
\smallskip

\noindent
{\bf Conjecture.}{\it
~Let ${\frak c}=\{c(\nu)\}$ be a sequence (not necessary nonnegative). 
Define $D_{\frak c}(s)$ and $Z_{\frak c}(x)$ by $\eqref{000}$ and $\eqref{009}$, respectively. 
Suppose that 
\begin{enumerate}
\item[(s-1)] $\frak c$ satisfies the growth condition (c-1), 
\item[(s-2)] $D_{\frak c}(s)$ is continued meromorphically to $\Bbb C$,   
\item[(s-3)] $D_{\frak c}(s)$ has double pole at $s=1$,
\item[(s-4)] $D_{\frak c}(s)$ has no poles in $(1/2,1)$, 
\item[(s-5)] $D_{\frak c}(s)$ has pole of even order at $s=1/2$,
\item[(s-6)] $D_{\frak c}(s)$ has a (suitable) functional equation for $s$ to $1-s$,  
\item[(s-7)] $D_{\frak c}(s)$ has a Euler product which converges absolutely for $\Re(s)>1$.
\end{enumerate} 
Then $Z_{\frak c}(x)$ has a single sign for sufficiently small $x>0$. 
}
\smallskip

Here we do not suppose the nonnegativety of $\{c(\nu)\}$, 
because it is not obtained from relations $\eqref{908}$ and $\eqref{910}$ in general. 
The nonnegativety of coefficients $\{c(\nu)\}$ of $D(s)$ 
is obtained by a suitable bound condition for $\{\omega(p^{k})\}_{p,k}$ in $\eqref{909}$. 

Anyway, it seems that it is difficult to prove the conjecture mentioned above, 
even if we suppose the further condition that $\{c(\nu)\}$ is nonnegative. 
At least, the author has no idea about it. 
Probably, to obtain a progress for the single sign property of $Z(x)$, 
we should study $Z(x)$ attached to specific nice Dirichlet series. 
For instance, $Z(x)$ attached to an elliptic curve, $Z(x)$ attached to
\[
D_{1,k}(s) 
  = \left( \frac{\zeta(2s)\zeta(2s-1)}{\zeta(2s-1/2)^k} \right)^2, \quad 
D_{\chi,k}(s) 
  = \left( \frac{\zeta(2s)\zeta(2s-1)}{L(2s-1/2,\chi)^k} \right)^2,
\]
where $L(s,\chi)$ is the Dirichlet $L$-function associated with 
the primitive Dirichlet character $\chi$, 
or $Z(x)$ attached to an arithmetic scheme (\cite[section 4.3]{IGM}).  

\subsection{Partial Euler product}
In section 3, we state that a finite truncation of the series expansion of $Z_{\frak c}(x)$ 
has a single sign for small $x>0$. 
In this part, we remark on another finitization of $Z_{\frak c}(x)$ in the case of elliptic curve 
according to Goldfeld~\cite{MR679556} and Conrad~\cite{MR2124918}.  
Let $E/{\Q}$ be an elliptic curve having conductor $q_E$. Define
\[
L_T(E,s) = 
\prod_{{p|q_E}\atop{p \leq T}} (1- a(p) p^{-s})^{-1}
\prod_{{p\not\,|q_E}\atop{p \leq T}} (1 - a(p) p^{-s} + p^{1-2s})^{-1}. 
\]
This gives the nonnegative sequence ${\frak c}=\{c_{E,T}(\nu)\}$ by 
\[
D_{E,T}(s) = \sum_{\nu} c_{E,T}(\nu) \, \nu^{-s} = \left( \frac{\zeta(2s)\zeta(2s-1)}{q_E^s L_T(E,2s)} \right)^2. 
\]
In this case, $c_{E,T}(\nu)\not=0$ for infinitely many $\nu \geq 1$. 
Further we find that for any fixed $T \geq 2$, 
there exists a negative constant $C_1(T)$ such that 
\begin{equation}\label{901}
Z_{\frak c}(x) = C_1(T) \, x \log x + O_T(x) \quad \text{as} \quad x \to 0^+.
\end{equation}
In fact the integral formula
\[
Z_{\frak c}(x) 
 = \frac{x^2}{2\pi i} \int_{(2+\delta)} 2^{-3}(2 \pi)^{s/2} \, (s-2)^4 \, \widehat{\zeta}(s/2)^2 D_{E,T}(s/2) \, x^{-s} ds
\]
and the residue theorem give
\[
\aligned
Z_{\frak c}(x) & = 
C_{1}(T) \, x \log x + \tilde{C}_{1}(T) \, x 
+ C_{0}(T) \, x^2 \log x + \tilde{C}_{0}(T) \, x^2 
\\
& \quad + \frac{x^2}{2\pi i}
\int_{(-\delta)}
2^{-3} (2 \pi)^{s/2} \, (s-2)^4 \, \widehat{\zeta}(s/2)^2 D_{E,T}(s/2) \, x^{-s} ds \\
& = 
C_{1}(T) \, x \log x + \tilde{C}_{1}(T) \, x 
+ C_{0}(T) \, x^2 \log x + \tilde{C}_{0}(T) \, x^2 + O_T(x^{2+\delta}),
\endaligned
\]
where 
\begin{align}\label{902}
C_{1}(T)=
- q_E^{-1} \frac{\Gamma(1/4)^2}{16 \sqrt{2} } \,  \,
\frac{ \zeta(0)^2 \zeta(1/2)^2}{ L_T(E,1)^2 }.
\end{align}
Hence, the expected equality obtained by the limit $T \to +\infty$ of $\eqref{610}$, 
$\eqref{901}$ and $\eqref{902}$ suggest that 
\begin{equation}\label{903}
L_T(E,1) \sim C \, (\log T)^{-r} \quad \text{as} \quad T \to \infty
\end{equation}
for some $0 \not= C \in {\R}$ and $r \geq 0$. 
Goldfeld~\cite{MR679556} proved that if $\eqref{903}$ holds, 
then $L(E,s)$ satisfies the Riemann hypothesis, the BSD-conjecture with $r = {\rm ord}_{s=1}L(E,s)$ and 
\[
C = \frac{L^{(r)}(E,1)}{r!} \cdot \frac{1}{\sqrt{2}\,e^{r\gamma}},
\]
where $\gamma=0.577215+$ is Euler's constant. 

Thus a detail study for the constant $C_1(T)$ and the error term $O_T(x)$ in $\eqref{901}$ 
may be available for the single sign property of $Z(x)$. 

\bibliographystyle{amsplain}
\bibliography{biblio}

\bigskip

\noindent
Masatoshi Suzuki\\
Department of Mathematics \\ 
Rikkyo University \\
Nishi-Ikebukuro, Toshima-ku \\
Tokyo 171-8501, \\
Japan \\
\texttt{suzuki@@rkmath.rikkyo.ac.jp}

\end{document}